\documentclass{arxivpaper}
\usepackage{manuscript}

\begin{document}
\maketitle
\thispagestyle{empty}
\begin{abstract}
This paper studies the shared-feedback design problem for finite-horizon averaged-cost control under model uncertainty represented by a joint law over the model parameter and initial state.  One feedback map is shared across all realizations and uses only time and the current state as online inputs.  We develop two methods: the nonlinear adjoint method combines forward closed-loop solves, backward adjoints, and joint-law policy gradients, while the linear moment method gives exact cost and gradient formulas for a homogeneous linear--quadratic design model by propagating conditional second moments.  Numerical studies compare averaged and nominal designs for time-dependent neural feedback and static linear feedback in nonlinear Van der Pol and Stuart--Landau systems; the static gains are designed from oscillator linearizations and evaluated on the original nonlinear systems.  Averaged design has lower cost in most comparisons.
\end{abstract}


\section{Introduction}

Feedback controllers intended for repeated deployment must operate despite variations in physical parameters, initial conditions, and operating environments.  Manufacturing tolerances alter physical dynamics, while simulation-to-hardware transfer introduces additional model mismatch \citep{ito2019manufacturing,peng2018sim}.  When the realized model parameter is unknown or unavailable online, designing and maintaining a separate controller for every realization is impractical.  This motivates a single feedback law that uses the observed state while balancing performance across possible models and initial states.

Several lines of research have addressed control under model and initial-state uncertainty.  They differ in what they optimize and in what the controller may use online.  One line chooses a single time-dependent control for all possible dynamics and evaluates its performance across them \citep{li2011ensemble,ross2015riemann,bettiol2019necessary,bettiol2024necessary,pesare2021convergence,aronna2026dynamic}.  Finite ensembles and sampled models make this problem tractable \citep{phelps2016optimal,lambrianides2020scalable,scagliotti2023optimal,aronna2025average,melnikov2026convergence}.  A second line optimizes a feedback law.  Some methods solve a Riccati problem for a finite collection of models and then restrict the resulting controller so that it acts on one system state at a time \citep{guth2024tracking,guth2025ensemble}.  Other methods allow the feedback law to use the realized model parameter \citep{verma2025neural}.  Methods that do not use the model parameter online rely only on the observed state.  For fixed dynamics, such feedback has been designed over uncertain initial states \citep{kunisch2024optimal}.  Methods that optimize one state-feedback law shared across model realizations have been developed for uncertain model parameters \citep{ito2024stochastic,ito2026stabilizing,fujinami2025policy}.

We study a finite-horizon problem in which the model parameter and initial state follow a joint law, while the model parameter is unavailable to the controller online.  Each realization is $y=(\omega,\xi)$, and a single feedback map $\pi_\theta(t,x)$ uses only time and the current state.  The objective averages the running and terminal costs over the joint law.  We call this the shared-feedback design problem and derive first-order conditions for two cases: nonlinear feedback with finitely many parameters and static linear feedback for homogeneous linear--quadratic systems.

The paper makes three contributions.
\begin{enumerate}[label=\textup{(\roman*)}]
\item We formulate the shared-feedback design problem under a joint law of model parameters and initial states.  Closely related studies optimize one state-feedback law shared across stochastic model realizations \citep{ito2024stochastic,ito2026stabilizing,fujinami2025policy}.  Their objectives are posed on an infinite horizon, and their feedback classes are time independent.  Our formulation is finite-horizon, allows $\pi_\theta(t,x)$, and uses a joint law that may couple $\omega$ and $\xi$.
\item For nonlinear feedback, we prove that the closed-loop state depends differentiably on the policy parameters and justify differentiating the objective after averaging over the joint law under explicit integrable bounds.  Building on adjoint methods for continuous-time policy optimization \citep{ainsworth2021faster,sandoval2023neural,kunisch2024optimal}, we obtain the gradient and first-order stationarity formulas used by the nonlinear adjoint method.
\item For homogeneous linear--quadratic systems with static feedback, we adapt classical Lyapunov and Gramian identities \citep{levine1970static,fatkhullin2021optimizing} to shared-gain design on a finite horizon, allowing time-dependent system and cost matrices, terminal costs, and model-dependent initial second moments.  The resulting linear moment method replaces initial-state particle propagation by fixed-dimensional moment equations after those moments have been formed.
\end{enumerate}
The remainder of the paper is organized as follows.  Section~\ref{sec:problem} formulates the shared-feedback design problem.  Section~\ref{sec:methods} develops the nonlinear adjoint method and the linear moment method.  Section~\ref{sec:numerical} applies both methods to Van der Pol and Stuart--Landau oscillators; under common evaluation laws, the controllers selected by averaged design have lower costs than those selected by nominal design in most comparisons.  Section~\ref{sec:conclusions} summarizes the results.

\paragraph{Notation.}
For a positive integer $d$, $\R^d$ is the $d$-dimensional Euclidean space and $I_d$ is its identity matrix.  For a matrix $A$, $A^\top$ and $\tr(A)$ denote its transpose and trace.  We write $\Sym(d)$ for the symmetric $d\times d$ matrices, and $A\succeq0$ and $A\succ0$ for positive semidefiniteness and positive definiteness.  The symbol $\norm{\cdot}$ denotes the Euclidean norm for vectors and the induced operator norm for matrices, while $\norm{v}_{\max}:=\max_i|v_i|$ and $\norm{A}_{\max}:=\max_{i,j}|A_{ij}|$; $\ip{\cdot}{\cdot}$ is the Euclidean inner product, $\ip{A}{B}_F:=\tr(A^\top B)$ is the Frobenius inner product, and $\norm{A}_F:=\sqrt{\ip{A}{A}_F}$.  For a topological space $E$, $\mathcal B_E$ is its Borel sigma-algebra, $\cP(E)$ is the set of Borel probability measures on $E$, and $\delta_z$ is the Dirac measure at $z$.  The symbols $D_xh$ and $\nabla_xh$ denote the derivative of a vector-valued map and the column gradient of a scalar-valued map, respectively; a dot denotes a time derivative, and $\interior C$ is the interior of $C$.

\section{Problem formulation}\label{sec:problem}

Let $S<T$ and let $[S,T]\subset\R$ be the corresponding compact interval, let $\Omega$ be a Polish space equipped with its Borel sigma-algebra $\mathcal B_\Omega$, and set $\Y:=\Omega\times\R^n$.  The joint realization is $y=(\omega,\xi)\in\Y$, with Borel probability law $\nu\in\cP(\Y)$.  Dependence between the model parameter $\omega$ and initial state $\xi$ is allowed.  Write $\rho$ for the $\Omega$-marginal of $\nu$.  For $t\in[S,T]$, $x\in\R^n$, $u\in\R^m$, and $\omega\in\Omega$, let $f(t,x,u,\omega)\in\R^n$ denote the vector field.  A shared feedback policy is a map $\pi_\theta:[S,T]\times\R^n\to\R^m$, parameterized by $\theta\in\Theta\subset\R^q$.  The policy uses time and the current state as its online inputs.  For each $y=(\omega,\xi)$, the corresponding closed-loop state and control satisfy
\begin{equation}\label{eq:closed-loop-dynamics}
\begin{aligned}
 \dot x^\theta(t,y)
 &=f\bigl(t,x^\theta(t,y),u^\theta(t,y),\omega\bigr),
 &x^\theta(S,y)&=\xi,\\
 u^\theta(t,y)&=\pi_\theta\bigl(t,x^\theta(t,y)\bigr).
\end{aligned}
\end{equation}
Realized control signals can differ across joint realizations because their states differ, allowing dependence in $\nu$ to affect the feedback response through the state.

Let $L(t,x,u,\omega)\in\R$ and $g(x,\omega)\in\R$ denote the running and terminal costs; like $f$, both may depend on the model parameter $\omega$.  Define the realization and averaged costs by
\begin{equation*}
 J_y(\theta):=g\bigl(x^\theta(T,y),\omega\bigr)
 +\int_S^T L\bigl(t,x^\theta(t,y),u^\theta(t,y),\omega\bigr)\,dt,
 \qquad \cJ_\nu(\theta):=\int_\Y J_y(\theta)\,\nu(dy).
\end{equation*}
A parameter $\theta$ is called admissible if the closed-loop problem \eqref{eq:closed-loop-dynamics} has a unique solution on $[S,T]$ for $\nu$-almost every $y$, the solution has a jointly measurable version $(t,y)\mapsto x^\theta(t,y)$, and $J_y(\theta)$ is measurable with integrable positive and negative parts.  The corresponding $\pi_\theta$ is an admissible shared policy.  Given a prescribed feasible set $\Theta\subset\R^q$ of admissible parameters, the shared-feedback design problem is
\begin{equation}\label{eq:shared-feedback-problem}
 \min_{\theta\in\Theta}\ \cJ_\nu(\theta).
\end{equation}
A fixed initial state $\xi_0$ is represented by $\nu(d\omega,d\xi)=\rho(d\omega)\delta_{\xi_0}(d\xi)$, while for an initial-state law $\mu\in\cP(\R^n)$, independent model and initial-state uncertainty has the product form $\nu(d\omega,d\xi)=\rho(d\omega)\mu(d\xi)$.  A nonproduct joint law represents dependence between the two.

\begin{remark}[Relation to existing uncertainty formulations]
For a finite uncertainty ensemble, average and ensemble open-loop formulations impose $u_i(t)=u(t)$, so that one signal is shared by all realizations \citep{phelps2016optimal,lambrianides2020scalable,bettiol2019necessary,bettiol2024necessary,scagliotti2023optimal,aronna2025average,aronna2026dynamic,melnikov2026convergence,pesare2021convergence}.  Ensemble Riccati constructions first derive a centralized feedback $u(t)=\kappa_\Sigma(t,X(t))$, with $X=(x_1^\top,\ldots,x_N^\top)^\top$, for an auxiliary extended system and then restrict it to a local feedback $u_i(t)=K_\Sigma(t,x_i(t))$ \citep{guth2024tracking,guth2025ensemble}.  Parameter-indexed feedback instead has the form $u_i(t)=\pi(t,x_i(t),\omega_i)$ \citep{verma2025neural}.  Problem~\eqref{eq:shared-feedback-problem} imposes $u_i(t)=\pi_\theta(t,x_i(t))$: every realization uses the same map, observes only its own state, and does not use the model parameter online.  The same information structure appears in feedback design over initial-state ensembles for fixed dynamics \citep{kunisch2024optimal} and over stochastic model parameters \citep{ito2024stochastic,fujinami2025policy}.  The work of \citet{ito2026stabilizing} also averages over random initial states, but through separate model and initial-state laws.  Here a single joint law covers both and may include dependence.  These model-uncertainty studies use structured cost or value representations for infinite-horizon feedback, whereas our finite-horizon formulation uses trajectory adjoints without a value approximation for feedback design.
\end{remark}

\section{Methods for the shared-feedback design problem}\label{sec:methods}

This section derives two methods from first-order necessary conditions for the shared-feedback design problem~\eqref{eq:shared-feedback-problem}.  The nonlinear adjoint method solves a state--adjoint pair for each node of a finite joint design law, whereas the linear moment method retains each model node and represents its conditional initial-state law by a second moment.  All proofs are given in Appendix~\ref{sec:proofs-appendix}.

\subsection{Nonlinear adjoint method}

The nonlinear adjoint method obtains the objective gradient by solving the closed-loop state equation forward and the adjoint equation backward for each realization and then averaging over the general joint law $\nu$.  This forward--backward construction follows continuous-time adjoint policy optimization \citep{ainsworth2021faster,sandoval2023neural,kunisch2024optimal}; here we formulate it under a joint model--initial-state law.  The following regularity conditions establish the state sensitivity used in the adjoint calculation.

The first-order formulas are local in the policy parameter.  Realization-dependent tubes and integrable envelopes provide the differentiability and domination needed for joint averaging.  The derivation uses standard results for parameter-dependent integral equations and Gronwall estimates \citep{coddington1955theory}, together with dominated convergence and Fubini--Tonelli interchange \citep{folland1999real}.

\begin{assumption}[Local regularity and integrable domination]\label{ass:regularity}
Fix an admissible reference parameter $\bar\theta$.  There is an open neighborhood $U\subset\R^q$ of $\bar\theta$ with the following properties.
\begin{enumerate}[label=\textup{(A\arabic*)}]
\item The maps $f$ and $L$ are jointly measurable in $(t,x,u,\omega)$ and continuously differentiable in $(x,u)$ for each $(t,\omega)$, with jointly measurable first derivatives.  The map $g$ is jointly measurable in $(x,\omega)$ and continuously differentiable in $x$, with jointly measurable gradient.  The policy is jointly measurable in $(t,x,\theta)$ and continuously differentiable in $(x,\theta)$ for each $t$, with jointly measurable derivatives $D_x\pi_\theta$ and $D_\theta\pi_\theta$.
\item Every $\theta\in U$ is admissible.  Moreover, there is a measurable set $\Y_0\subset\Y$ with $\nu(\Y_0)=1$ such that, for every $y\in\Y_0$ and every $\theta\in U$, the closed-loop solution exists uniquely on $[S,T]$.  The trajectories, the line segments between trajectories used in parameter difference quotients, and the corresponding policy values lie in measurable realization-dependent tubes on which the mean-value formulas below are valid.
\item There are nonnegative measurable envelopes $a,c,b,d$ on $[S,T]\times\Y$ and $e$ on $\Y$.  For every $y\in\Y_0$, the four time-dependent envelopes are integrable on $[S,T]$ and $e(y)<\infty$.  On the tubes and for $\theta\in U$, the bounds $\norm{D_xf+D_ufD_x\pi_\theta}\le a$ and $\norm{D_ufD_\theta\pi_\theta}\le c$ hold for the dynamics, while $\norm{\nabla_xL+(D_x\pi_\theta)^\top\nabla_uL}\le b$ and $\norm{(D_\theta\pi_\theta)^\top\nabla_uL}\le d$ hold for the running cost.  Here derivatives of $f$ and $L$ are evaluated at $(t,x,\pi_\theta(t,x),\omega)$; in addition, $\norm{\nabla_xg}\le e$, with the envelopes evaluated at $(t,y)$ where applicable.
\item The following function of $y$ is integrable with respect to $\nu$:
\begin{equation}\label{eq:dominating-envelope}
 \left(e(y)+\int_S^T b(t,y)\,dt\right)
 \exp\!\left(\int_S^T a(t,y)\,dt\right)
 \int_S^T c(t,y)\,dt
 +\int_S^T d(t,y)\,dt.
\end{equation}
\end{enumerate}
\end{assumption}

\begin{remark}[Compact-support sufficient conditions]
If \textup{(A1)--(A2)} hold, $\nu$ has compact support, and one deterministic compact set contains all nearby trajectories and policy values, uniform derivative bounds by time-integrable functions imply \textup{(A3)} and the $\nu$-integrability of \eqref{eq:dominating-envelope}.  More generally, the tubes in \textup{(A2)} may depend measurably on $y$; unbounded laws are covered when global existence holds, the bounds in \textup{(A3)} are valid, and \eqref{eq:dominating-envelope} is $\nu$-integrable.
\end{remark}

For an adjoint variable $p\in\R^n$, define the Hamiltonian $H(t,x,u,p,\omega):=p^\top f(t,x,u,\omega)-L(t,x,u,\omega)$, whose control gradient is $\nabla_uH=D_uf^\top p-\nabla_uL$.  Unless stated otherwise, derivatives in this section are evaluated along the reference trajectory $x^{\bar\theta}(t,y)$ and the associated control $u^{\bar\theta}(t,y)$.

\begin{theorem}[Adjoint-based policy gradient]\label{thm:policy-gradient}
Suppose Assumption~\ref{ass:regularity} holds.  For $\nu$-almost every $y$, there is a unique absolutely continuous adjoint $p_y:[S,T]\to\R^n$ satisfying
\begin{equation}\label{eq:closed-loop-adjoint}
 \begin{cases}
 -\dot p_y
 =\bigl(D_xf+D_ufD_x\pi_{\bar\theta}\bigr)^\top p_y
  -\nabla_xL-(D_x\pi_{\bar\theta})^\top\nabla_uL,\\[1mm]
 -p_y(T)=\nabla_xg\bigl(x^{\bar\theta}(T,y),\omega\bigr).
 \end{cases}
\end{equation}
There is a jointly measurable version of the family $(t,y)\mapsto p_y(t)$ for which the following pairings are integrable.  The averaged objective is differentiable at $\bar\theta$, and
\begin{equation}\label{eq:policy-gradient}
 \nabla_\theta\cJ_\nu(\bar\theta)
 =\int_\Y\int_S^T(D_\theta\pi_{\bar\theta})^\top
 \bigl(\nabla_uL-D_uf^\top p_y\bigr)\,dt\,\nu(dy).
\end{equation}
\end{theorem}

For each realization $y=(\omega,\xi)$, write $(\nabla_u H)^{\bar\theta}(t,y):=\nabla_u H(t,x^{\bar\theta}(t,y),u^{\bar\theta}(t,y),p_y(t),\omega)$ for the control gradient evaluated along the corresponding state, control, and adjoint.  In Hamiltonian form, \eqref{eq:policy-gradient} is $\nabla_\theta\cJ_\nu(\bar\theta)=-\int_\Y\int_S^T(D_\theta\pi_{\bar\theta})^\top(\nabla_u H)^{\bar\theta}(t,y)\,dt\,\nu(dy)$.  At a local minimizer $\theta^*$, let $p_y^*$ be the corresponding adjoint and write $(\nabla_u H)^*(t,y):=\nabla_u H(t,x^{\theta^*}(t,y),u^{\theta^*}(t,y),p_y^*(t),\omega)$.  The resulting gradient gives the local stationarity condition for the selected policy class.
\begin{corollary}[Policy-parameter stationarity]\label{cor:projected-stationarity}
Suppose Assumption~\ref{ass:regularity} holds at a local minimizer $\theta^*$, so that Theorem~\ref{thm:policy-gradient} applies.  If $\theta^*\in\interior\Theta$, then
\begin{equation}\label{eq:param-stationarity}
 \int_\Y\int_S^T
 (D_\theta\pi_{\theta^*})^\top(\nabla_u H)^*(t,y)\,dt\,\nu(dy)=0.
\end{equation}
If instead $\Theta$ is closed and convex and $\theta^*$ is a local minimizer relative to $\Theta$, then
\begin{equation}\label{eq:param-VI}
 \ip{\nabla_\theta\cJ_\nu(\theta^*)}{\theta-\theta^*}\ge0
 \qquad(\theta\in\Theta).
\end{equation}
\end{corollary}

For numerical optimization, given quadrature or empirical nodes $y_j=(\omega_j,\xi_j)\in\Y$ and positive weights $w_j$ summing to one, use
\begin{equation}\label{eq:design-rule}
 \nu=\sum_{j=1}^{M}w_j\delta_{y_j},
 \qquad
 \cJ_\nu(\theta)=\sum_{j=1}^{M}w_j J_{y_j}(\theta).
\end{equation}
For this finite rule, items~\textup{(A1)--(A3)} of Assumption~\ref{ass:regularity} are imposed at the design nodes; finite support then makes item~\textup{(A4)} automatic, and \eqref{eq:policy-gradient} becomes a weighted sum.  We use the resulting objective and gradient, with $\bar\theta=\theta$ at each iterate, for policy updates and assess first-order convergence through
\begin{equation}\label{eq:policy-stationarity-residual}
 r_\Theta(\theta)
 :=\norm{\theta-\Pi_\Theta\bigl(\theta-\nabla_\theta\cJ_\nu(\theta)\bigr)}_{\max},
\end{equation}
where $\Pi_\Theta$ is the Euclidean projection onto the nonempty closed convex parameter set $\Theta$.  The finite-law variational inequality corresponding to \eqref{eq:param-VI} is equivalent to $r_\Theta(\theta^*)=0$, which gives the stopping test in Algorithm~\ref{alg:nonlinear-design}.

\begin{algorithm}[H]
\algocfconts{alg:nonlinear-design}{%
\caption{Nonlinear adjoint method for the shared-feedback design problem}%
}{%
\begin{algorithmic}[1]
\REQUIRE Fixed finite design law $\nu$; nonempty closed convex parameter set $\Theta$; initial parameter $\theta_0\in\Theta$; stationarity tolerance $\varepsilon_\theta>0$
\ENSURE Upon termination, a stationarity-certified policy parameter $\widehat\theta$ with $r_\Theta(\widehat\theta)\le\varepsilon_\theta$
\STATE Set $\theta\gets\theta_0$ and $r_\Theta(\theta)\gets+\infty$
\WHILE{$r_\Theta(\theta)>\varepsilon_\theta$}
  \FOR{each design realization $y_j$}
    \STATE Solve the closed-loop state equation~\eqref{eq:closed-loop-dynamics} forward
    \STATE Solve the adjoint equation~\eqref{eq:closed-loop-adjoint} backward along that trajectory
    \STATE Accumulate the cost and policy-gradient contribution~\eqref{eq:policy-gradient}
  \ENDFOR
  \STATE Evaluate $\cJ_\nu(\theta)$, $\nabla_\theta\cJ_\nu(\theta)$, and the residual~\eqref{eq:policy-stationarity-residual}
  \IF{$r_\Theta(\theta)>\varepsilon_\theta$}
    \STATE Apply one constrained optimization step to update $\theta\in\Theta$
  \ENDIF
\ENDWHILE
\STATE Return $\widehat\theta\gets\theta$
\end{algorithmic}
}
\end{algorithm}

\subsection{Linear moment method}

The linear moment method specializes the shared-feedback design problem~\eqref{eq:shared-feedback-problem} to homogeneous linear dynamics, quadratic costs, and static linear feedback.  Related work on linear systems incorporates model distributions through stochastic-parameter control design and averaged or domain-randomized {LQR} formulations \citep{fisher2009linear,pesare2021convergence,ito2024stochastic,fujinami2025policy}.  Shared-gain {LQR} with task-dependent initial covariances is studied in discrete time by \citet{stamouli2025multitask}, while \citet{cui2026probabilistic} develop a continuous-time gradient method using polynomial-chaos approximations.  Both consider an infinite horizon.  Here, we treat finite-horizon systems with time-dependent data and terminal costs, retaining a finite model marginal and integrating conditional initial laws exactly through their second moments.  The derivations use standard state-transition-matrix theory \citep{sontag1998mathematical} and finite-dimensional matrix variations.

Assume that the $\Omega$-marginal introduced in Section~\ref{sec:problem} is $\rho=\sum_{i=1}^N\alpha_i\delta_{\omega_i}$, where $\alpha_i>0$ and $\sum_{i=1}^N\alpha_i=1$, and write the joint law as $\nu(d\omega,d\xi)=\sum_{i=1}^N\alpha_i\delta_{\omega_i}(d\omega)\mu_i(d\xi)$, where $\mu_i\in\cP(\R^n)$ is the conditional initial law given $\omega_i$.  Dependence between the model and initial state is represented by the variation of $\mu_i$ with $i$.  For $y=(\omega_i,\xi)$, specialize the dynamics and costs of Section~\ref{sec:problem} by setting $f(t,x,u,\omega_i)=A_i(t)x+B_i(t)u$, $L(t,x,u,\omega_i)=\tfrac12(x^\top Q_i(t)x+u^\top R_i(t)u)$, and $g(x,\omega_i)=\tfrac12x^\top G_ix$, where $Q_i(t)\succeq0$, $R_i(t)\succeq rI_m$ for some $r>0$, and $G_i\succeq0$.  Specialize the shared policy to $\pi_K(t,x)=-Kx$, where $K\in\R^{m\times n}$ plays the role of $\theta$ (equivalently, $\theta=\operatorname{vec}(K)$ and $q=mn$).  Write $A_{i,K}(t):=A_i(t)-B_i(t)K$ and let $x_i(\cdot,\xi)$ denote the corresponding closed-loop trajectory.
The following assumption gives concrete sufficient conditions for Assumption~\ref{ass:regularity} in this linear--quadratic setting.
\begin{assumption}[Linear--quadratic integrability]\label{ass:lq-integrability}
For every $i$, the functions $A_i,B_i,Q_i$, and $R_i$ are essentially bounded on $[S,T]$, and $\int_{\R^n}\norm{\xi}^2\,\mu_i(d\xi)<\infty$.
\end{assumption}

\begin{proposition}[Conditional moment propagation and cost]\label{prop:conditional-moments}
Under Assumption~\ref{ass:lq-integrability}, define
\begin{equation}\label{eq:state-moment-definition}
 M_i^0:=\int_{\R^n}\xi\xi^\top\,\mu_i(d\xi),
 \qquad
 M_i(t):=\int_{\R^n}x_i(t,\xi)x_i(t,\xi)^\top\,\mu_i(d\xi).
\end{equation}
Then $M_i$ is absolutely continuous and, for almost every $t\in[S,T]$,
\begin{equation}\label{eq:moment-forward}
 \dot M_i=A_{i,K}M_i+M_iA_{i,K}^\top,
 \qquad M_i(S)=M_i^0,
\end{equation}
and the exact averaged cost is
\begin{equation*}
 \cJ_\nu(K)=\frac12\sum_{i=1}^N\alpha_i\left[
 \tr(G_iM_i(T))
 +\int_S^T\tr\bigl((Q_i+K^\top R_iK)M_i\bigr)\,dt
 \right].
\end{equation*}
\end{proposition}

\begin{proposition}[Lyapunov evaluation and value]\label{prop:policy-evaluation}
Under Assumption~\ref{ass:lq-integrability}, for a prescribed gain $K$, let $P_i:[S,T]\to\Sym(n)$ be the unique absolutely continuous solution of
\begin{equation}\label{eq:P-Lyapunov}
 \begin{cases}
 -\dot P_i=A_{i,K}^\top P_i+P_iA_{i,K}+Q_i+K^\top R_iK \quad\text{a.e. on }[S,T],\\
 P_i(T)=G_i.
 \end{cases}
\end{equation}
Then $p_i(t,\xi):=-P_i(t)x_i(t,\xi)$ solves the closed-loop adjoint equation~\eqref{eq:closed-loop-adjoint} for $y=(\omega_i,\xi)$, and
\begin{equation}\label{eq:trace-value-identity}
 \cJ_\nu(K)=\frac12\sum_{i=1}^N\alpha_i\tr\bigl(P_i(S)M_i^0\bigr).
\end{equation}
\end{proposition}

For fixed $K$, \eqref{eq:P-Lyapunov} is the linear Lyapunov equation for policy evaluation.

\begin{theorem}[Moment-based shared-gain gradient]\label{thm:shared-gain-derivative}
Under Assumption~\ref{ass:lq-integrability}, for a prescribed constant gain $K\in\R^{m\times n}$, let $M_i$ and $P_i$ solve \eqref{eq:moment-forward} and \eqref{eq:P-Lyapunov}.  The objective is differentiable with respect to $K$, and
\begin{equation}\label{eq:gain-derivative}
 \nabla_K\cJ_\nu(K)
 =\int_S^T\sum_{i=1}^N\alpha_i
 (R_iK-B_i^\top P_i)M_i\,dt.
\end{equation}
\end{theorem}

The Lyapunov/moment gradient structure is classical in static-feedback optimization \citep{levine1970static,fatkhullin2021optimizing}.  Here it is aggregated over the model nodes and their conditional initial moments.

\begin{corollary}[Shared-gain stationarity]\label{cor:shared-gain-stationarity}
Under Assumption~\ref{ass:lq-integrability}, let $P_i^*$ and $M_i^*$ denote the solutions evaluated at $K=K^*$.  Every local minimizer $K^*\in\R^{m\times n}$ satisfies
\begin{equation}\label{eq:shared-gain-stationarity}
 \int_S^T\sum_{i=1}^N\alpha_i
 \bigl(R_i(t)K^*-B_i(t)^\top P_i^*(t)\bigr)M_i^*(t)\,dt=0.
\end{equation}
\end{corollary}
Define the first-order residual
\begin{equation}\label{eq:gain-stationarity-residual}
 r(K):=\norm{\nabla_K\cJ_\nu(K)}_{\max}.
\end{equation}

\begin{corollary}[Shared-gain fixed-point identity under excitation]\label{cor:shared-gain-fixed-point}
Under Assumption~\ref{ass:lq-integrability}, assume in addition that $R_i(t)=R\succ0$ is common and constant.  Define
\begin{equation*}
 \mathsf S(K):=\sum_{i=1}^N\alpha_i\int_S^T M_i(t)\,dt,
 \qquad
 \mathsf B(K):=\sum_{i=1}^N\alpha_i\int_S^T B_i(t)^\top P_i(t)M_i(t)\,dt.
\end{equation*}
The matrix $\mathsf S(K)$ is positive definite if and only if
\begin{equation}\label{eq:excitation-condition}
 \sum_{i=1}^N\alpha_i\int_S^T
 \int_{\R^n}\bigl(v^\top x_i(t,\xi)\bigr)^2\,\mu_i(d\xi)dt>0
 \qquad\text{for every }v\ne0.
\end{equation}
If \eqref{eq:excitation-condition} holds, then $K$ is stationary if and only if
\begin{equation}\label{eq:shared-gain-fixed-point}
 K=R^{-1}\mathsf B(K)\mathsf S(K)^{-1}.
\end{equation}
The Anderson--Moore direction
\begin{equation}\label{eq:anderson-moore-direction}
 D(K):=R^{-1}\mathsf B(K)\mathsf S(K)^{-1}-K
 =-R^{-1}\nabla_K\cJ_\nu(K)\mathsf S(K)^{-1}
\end{equation}
is a strict descent direction whenever $\nabla_K\cJ_\nu(K)\ne0$.
\end{corollary}

A simple sufficient condition for \eqref{eq:excitation-condition} is $M_i^0\succ0$ for at least one index with positive weight.  Equation~\eqref{eq:shared-gain-fixed-point} is an implicit fixed-point identity because $P_i$ and $M_i$ depend on $K$.  Under the common-$R$ and excitation conditions of Corollary~\ref{cor:shared-gain-fixed-point}, Algorithm~\ref{alg:linear-design} uses the corresponding Anderson--Moore direction with Armijo backtracking, following the classical static-feedback iteration viewpoint \citep{makila1984gradient,toivonen1985algorithm}.

\begin{algorithm}[H]
\algocfconts{alg:linear-design}{%
\caption{Linear moment method for the shared-feedback design problem}%
}{%
\begin{algorithmic}[1]
\REQUIRE Model nodes and weights $(\omega_i,\alpha_i)$; conditional initial moments $M_i^0$; common constant $R_i(t)=R\succ0$; initial gain $K_0$; tolerance $\varepsilon_K>0$; excitation condition \eqref{eq:excitation-condition} at every iterate; Anderson--Moore update routine with Armijo backtracking
\ENSURE Upon termination, a gain $\widehat K$ for static feedback with first-order residual at most $\varepsilon_K$
\STATE Set $K\gets K_0$
\STATE Evaluate $K$ by solving \eqref{eq:moment-forward} and \eqref{eq:P-Lyapunov} for $i=1,\ldots,N$ and forming $\cJ_\nu(K)$, $\nabla_K\cJ_\nu(K)$, $r(K)$, $\mathsf S(K)$, and $\mathsf B(K)$
\WHILE{$r(K)>\varepsilon_K$}
  \STATE Apply one Anderson--Moore update based on \eqref{eq:anderson-moore-direction}, using Armijo backtracking, to obtain a new $K$
  \STATE Re-evaluate the updated gain as above
\ENDWHILE
\STATE Return $\widehat K\gets K$
\end{algorithmic}
}
\end{algorithm}

\begin{remark}[Moment closure and structured iteration]\label{rem:moment-computation}
Since the state is linear in $\xi$, the value and gain derivative depend on the initial law only through $M_i^0$, yielding the exact closure \eqref{eq:moment-forward}--\eqref{eq:P-Lyapunov}.  After these moments are formed, propagation uses $Nn^2$ variables rather than $NMn$ for particles, independently of $M$.  Given $P_i$ and $M_i$, \eqref{eq:anderson-moore-direction} computes the gain direction from $\mathsf S(K)$ and $\mathsf B(K)$ without a generic constrained optimizer or quasi-Newton model.  The nonlinear adjoint method uses realization-wise state and adjoint solves and a general policy-parameter update.
\end{remark}

\section{Numerical study}\label{sec:numerical}

We conduct numerical studies on two models and compare averaged design with nominal design.  Averaged design uses the prescribed model--initial-state law, whereas nominal design uses $\nu_{\mathrm{nom}}:=\delta_{(\bar\omega,\bar\xi)}$, where $(\bar\omega,\bar\xi)=(\mathbb E[\omega],\mathbb E[\xi])$.  For each model, we apply the nonlinear adjoint method to the nonlinear dynamics with time-dependent neural feedback and the linear moment method to the linearized dynamics with static feedback.  The computational environment and implementation details are reported in Appendix~\ref{sec:numerical-appendix}.

Both models use $[S,T]=[0,5]$ and scalar actuation in the second state equation.  The controlled Van der Pol dynamics are $\dot q=v$ and $\dot v=-q+\omega(1-q^2)v+u$ \citep{vanderpol1926relaxation}, with $\omega\sim\mathrm{Uniform}[-2,4]$, $\xi=(1,1)^\top+\zeta$, $\zeta\sim\mathrm{Uniform}([-0.2,0.2]^2)$, and the running cost is $L(x,u)=\tfrac12(x^\top x+0.2u^2)$.  The controlled Stuart--Landau dynamics are $\dot q=\omega q-v-(q^2+v^2)q$ and $\dot v=q+\omega v-(q^2+v^2)v+u$ \citep{landau1965turbulence,stuart1958nonlinear}, with $\omega\sim\mathrm{Uniform}[-2.1,0.8]$, $\xi=(0.5,0.5)^\top+\zeta$, $\zeta\sim\mathrm{Uniform}([-0.2,0.2]^2)$, and the running cost is $L(x,u)=\tfrac12(2x^\top x+u^2)$.  Both experiments use zero terminal cost, $g\equiv0$.  In both experiments, the model parameter and initial-state perturbation are independent.  For the comparisons using static linear feedback, the design models are linearized about the zero equilibrium.

\begin{figure}[H]
\centering
\includegraphics[width=0.85\linewidth]{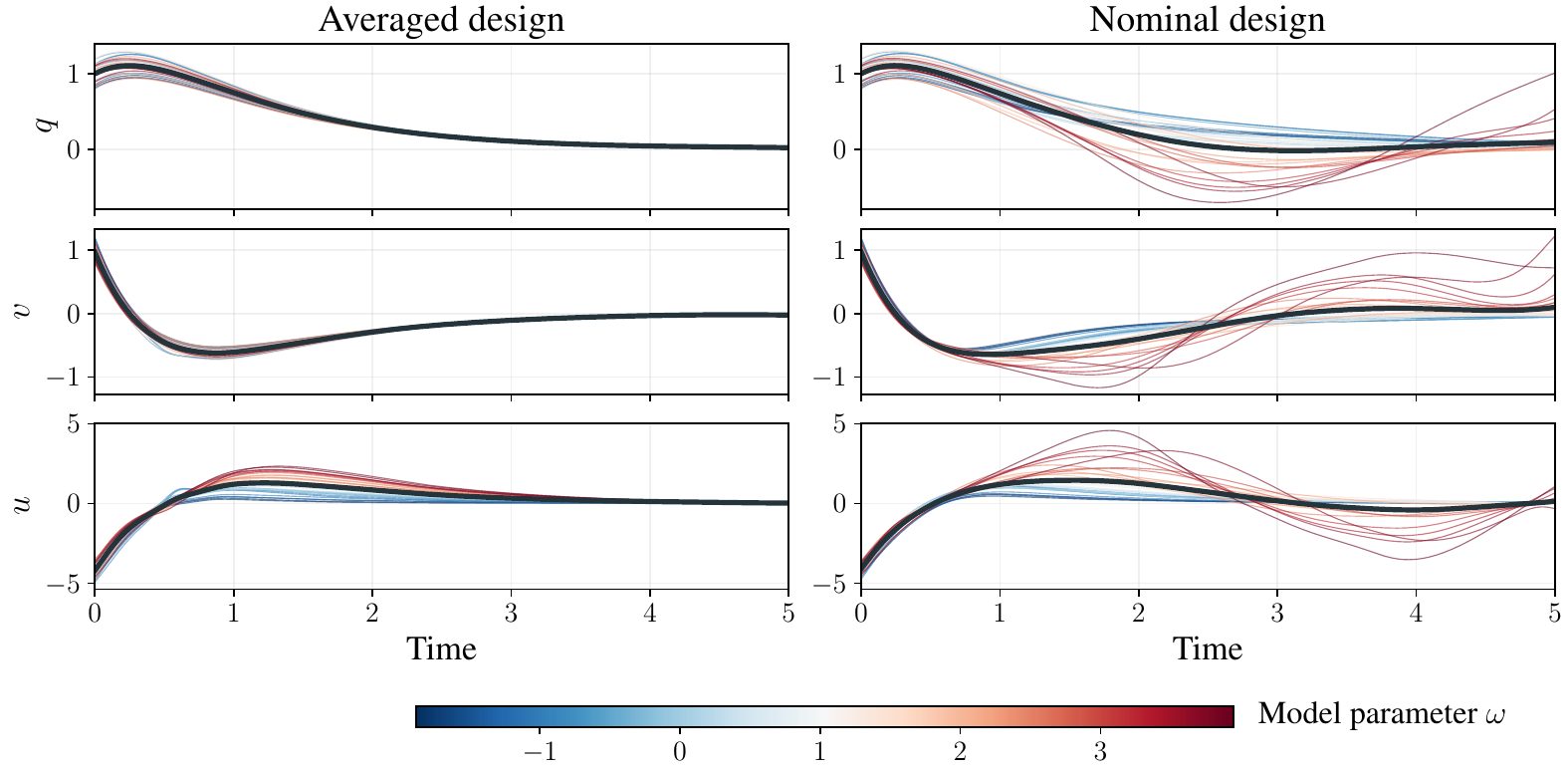}
\caption{Nonlinear Van der Pol trajectories under time-dependent neural feedback obtained by averaged design (left) and nominal design (right).  Both panels use the same 32 pairs of model parameters and initial states.  Rows show $q$, $v$, and $u$; color encodes the model parameter, and thick curves are equal-weight pointwise means.}
\label{fig:nonlinear-vdp-trajectories}
\end{figure}

Figure~\ref{fig:nonlinear-vdp-trajectories} compares nonlinear Van der Pol trajectories under time-dependent neural feedback for the same 32 pairs of model parameters and initial states.  Color identifies $\omega$; variation across curves also reflects $\xi$.  Averaged design yields a tighter state-and-control bundle and suppresses the parameter-dependent oscillations under nominal design.  Appendix~\ref{sec:numerical-appendix} reports the Stuart--Landau trajectories under time-dependent neural feedback and the two static-linear-feedback comparisons.

\begin{table}[t]
\centering
\caption{Finite-horizon averaged costs for the Van der Pol and Stuart--Landau systems.  Relative reduction (\%) is $100(J_{\mathrm{nom}}-J_{\mathrm{avg}})/J_{\mathrm{nom}}$.}
\label{tab:evaluation-costs}
\small
\begin{tabular}{@{}lllrrr@{}}
\toprule
System & Design dynamics & Feedback & $J_{\mathrm{avg}}$ & $J_{\mathrm{nom}}$ & Reduction\\
\midrule
\multirow{2}{*}{Nonlinear Van der Pol}
 & Nonlinear & Time-dependent neural & 1.379 & 1.765 & 21.9\%\\
 & Linearized & Static linear & 1.563 & 1.656 & 5.62\%\\
\multirow{2}{*}{Nonlinear Stuart--Landau}
 & Nonlinear & Time-dependent neural & 0.2845 & 0.3953 & 28.0\%\\
 & Linearized & Static linear & 0.2952 & 0.3878 & 23.9\%\\
\bottomrule
\end{tabular}
\end{table}

Table~\ref{tab:evaluation-costs} compares all four controllers for each nonlinear oscillator under a common evaluation law.  Averaged design has lower cost in both feedback classes, with the averaged neural policy attaining the lowest cost for each oscillator.  Appendix~\ref{sec:numerical-appendix} varies model-parameter and initial-state uncertainty separately.

\begin{figure}[t]
\centering
\includegraphics[width=0.85\linewidth]{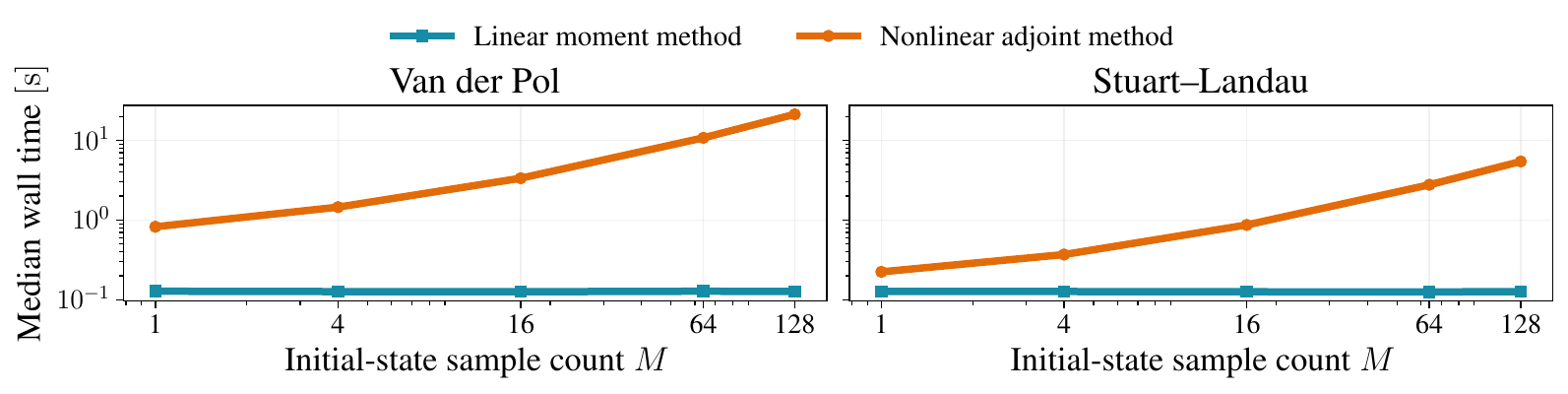}
\caption{Median wall time required to evaluate the design objective and its gradient once as the initial-state sample count increases.  Both methods use the same 32 model-parameter nodes and weights.  Blue denotes the linear moment method with static linear feedback, and orange denotes the nonlinear adjoint method with time-dependent neural feedback.}
\label{fig:runtime-scaling}
\end{figure}

Figure~\ref{fig:runtime-scaling} reports the time to evaluate the design objective and its gradient once using the same model-parameter nodes and weights.  On the nonlinear design model, the adjoint method propagates a state--adjoint pair for each $(\omega_i,\xi_j)$.  On the linearized design model, the moment method forms $\widehat M_i^0=M^{-1}\sum_{j=1}^M\xi_j\xi_j^\top$ and solves \eqref{eq:moment-forward} in a dimension independent of $M$.  Its runtime remains nearly constant over the tested sample counts, whereas the adjoint runtime increases with $M$.

\section{Conclusions}\label{sec:conclusions}

We studied shared-feedback design under a joint model--initial-state law, computing nonlinear policy gradients through trajectory adjoints.  For the linear--quadratic design model, the linear moment method computes exact costs and gradients for static linear feedback using conditional moment propagation and Lyapunov equations.  In the numerical study, averaged design has lower cost than nominal design in most comparisons across the two feedback classes and oscillator systems.  Exploring sufficient optimality conditions and studying convergence of the two design algorithms are left for future work.

\FloatBarrier
\subsection*{Generative AI Use}
Generative AI tools were used to assist with numerical computations and writing. The authors reviewed and verified all AI-assisted content and take responsibility for the final content of the paper.
\subsection*{Acknowledgements}
This project has received financial support from the CNRS through the MITI interdisciplinary programs.
D. Kalise is supported by the Air Force Office of Scientific Research under award number FA8655-26-1-B011.
M. G\'omez-Aedo was supported by ANID (Becas Chile, 2023, No. 72230073).

\appendix
\section{Proofs of the results}\label{sec:proofs-appendix}

This appendix contains the proofs of the results stated in the main text.  The state sensitivity used in the policy-gradient proof is stated here as an auxiliary lemma; all other notation, assumptions, and numbered statements refer to the main text.  For the proofs, abbreviate the closed-loop vector field and running cost by
\begin{equation*}
 F_\theta(t,x,\omega):=f(t,x,\pi_\theta(t,x),\omega),
 \qquad
 \ell_\theta(t,x,\omega):=L(t,x,\pi_\theta(t,x),\omega).
\end{equation*}

\subsection{Auxiliary state sensitivity for the policy-gradient theorem}

The proof of Theorem~\ref{thm:policy-gradient} uses the following auxiliary result.
\begin{lemma}[Closed-loop state sensitivity]\label{lem:state-sensitivity}
Under items~\textup{(A1)--(A3)} of Assumption~\ref{ass:regularity}, for every $y\in\Y_0$ the map $\theta\mapsto x^\theta(\cdot,y)$ from $U$ into $C([S,T];\R^n)$ is Fr\'echet differentiable at $\bar\theta$.  For a direction $\eta\in\R^q$, define
\begin{equation}\label{eq:state-variation-definition}
 z_y(t):=D_\theta x^{\bar\theta}(t,y)[\eta].
\end{equation}
Then $z_y$ is the unique solution of
\begin{equation}\label{eq:state-variation}
 \dot z_y(t)=A_y(t)z_y(t)+C_y(t)\eta,
 \qquad z_y(S)=0,
\end{equation}
where, along the reference trajectory, $A_y=D_xf+D_ufD_x\pi_{\bar\theta}\in\R^{n\times n}$ and $C_y=D_ufD_\theta\pi_{\bar\theta}\in\R^{n\times q}$.  Moreover,
\begin{equation}\label{eq:state-variation-bound}
 \sup_{t\in[S,T]}\norm{z_y(t)}
 \le
 \exp\!\left(\int_S^T a(t,y)\,dt\right)
 \int_S^T c(t,y)\,dt\,\norm{\eta}.
\end{equation}
\end{lemma}

\begin{resultproof}{Lemma~\ref{lem:state-sensitivity}}
Fix a realization $y=(\omega,\xi)\in\Y_0$ for which the bounds in item~\textup{(A3)} hold, and abbreviate
\begin{equation*}
 \bar x(t):=x^{\bar\theta}(t,y),
 \qquad
 x_h(t):=x^{\bar\theta+h}(t,y),
 \qquad
 \delta_h(t):=x_h(t)-\bar x(t).
\end{equation*}
Here $h\in\R^q$ is small enough that $\bar\theta+h\in U$.  Because the initial coordinate does not depend on the policy parameter,
$\delta_h(S)=0$.  Subtracting the two integral equations in
\eqref{eq:closed-loop-dynamics} gives
\begin{equation}\label{eq:state-difference-integral}
 \delta_h(t)
 =\int_S^t
 \left[F_{\bar\theta+h}(s,x_h(s),\omega)
       -F_{\bar\theta}(s,\bar x(s),\omega)\right]ds.
\end{equation}
Split the difference in the integrand into a state increment and a parameter increment.  The mean-value formula yields
\begin{align}
 F_{\bar\theta+h}(t,x_h,\omega)-F_{\bar\theta}(t,\bar x,\omega)
 &=A_h(t)\delta_h(t)+C_h(t)h,\label{eq:state-difference-mvf}\\
 A_h(t)
 &:=\int_0^1D_xF_{\bar\theta+h}
       \bigl(t,\bar x(t)+s\delta_h(t),\omega\bigr)\,ds,\nonumber\\
 C_h(t)
 &:=\int_0^1D_\theta F_{\bar\theta+sh}
       \bigl(t,\bar x(t),\omega\bigr)\,ds.\nonumber
\end{align}
All points in these integrals lie in the tube specified in Assumption~\ref{ass:regularity}.  Hence
$\lVert A_h(t)\rVert\le a(t,y)$ and
$\lVert C_h(t)\rVert\le c(t,y)$ for almost every $t$.  Equations
\eqref{eq:state-difference-integral}--\eqref{eq:state-difference-mvf} and Gronwall's inequality \citep{coddington1955theory} therefore imply
\begin{equation}\label{eq:state-lipschitz-bound}
 \sup_{t\in[S,T]}\lVert\delta_h(t)\rVert
 \le
 \exp\!\left(\int_S^T a(s,y)\,ds\right)
 \int_S^T c(s,y)\,ds\,\lVert h\rVert.
\end{equation}
In particular, $x_h\to\bar x$ uniformly as $h\to0$.

Set
\begin{equation*}
 A_y(t):=D_xF_{\bar\theta}(t,\bar x(t),\omega),
 \qquad
 C_y(t):=D_\theta F_{\bar\theta}(t,\bar x(t),\omega).
\end{equation*}
The chain rule for the closed-loop vector field gives
\begin{equation*}
 A_y=D_xf+D_ufD_x\pi_{\bar\theta},
 \qquad
 C_y=D_ufD_\theta\pi_{\bar\theta}.
\end{equation*}
For each $h$, let $z[h]$ denote the solution of
\begin{equation}\label{eq:linearized-state-h}
 \dot z[h]=A_yz[h]+C_yh,
 \qquad z[h](S)=0.
\end{equation}
Variation of constants, followed by Gronwall's inequality, gives
\begin{equation}\label{eq:linearized-state-bound-h}
 \sup_{t\in[S,T]}\lVert z[h](t)\rVert
 \le
 \exp\!\left(\int_S^T a(s,y)\,ds\right)
 \int_S^T c(s,y)\,ds\,\lVert h\rVert.
\end{equation}

It remains to verify that $z[h]$ is the derivative rather than only a formal linearization.  Define
$r_h:=\delta_h-z[h]$.  Subtracting \eqref{eq:linearized-state-h} from the integral form of
\eqref{eq:state-difference-mvf} gives
\begin{equation}\label{eq:state-remainder-integral}
 r_h(t)=\int_S^t\left[
 A_h(s)r_h(s)
 +(A_h(s)-A_y(s))z[h](s)
 +(C_h(s)-C_y(s))h
 \right]ds.
\end{equation}
The continuity assumptions, the uniform convergence $x_h\to\bar x$, and the tube condition imply
$A_h(t)\to A_y(t)$ and $C_h(t)\to C_y(t)$ for almost every $t$.  Moreover,
$\lVert A_h-A_y\rVert\le2a(\cdot,y)$ and
$\lVert C_h-C_y\rVert\le2c(\cdot,y)$.  Since $a(\cdot,y)$ and $c(\cdot,y)$ are integrable, dominated convergence gives
\begin{equation}\label{eq:coefficient-L1-convergence}
 \int_S^T\lVert A_h-A_y\rVert\,dt\longrightarrow0,
 \qquad
 \int_S^T\lVert C_h-C_y\rVert\,dt\longrightarrow0.
\end{equation}
Using \eqref{eq:linearized-state-bound-h} and \eqref{eq:coefficient-L1-convergence} in \eqref{eq:state-remainder-integral}, and then applying Gronwall once more, yields
\begin{equation}\label{eq:state-Frechet-remainder}
 \frac{\sup_{t\in[S,T]}\lVert r_h(t)\rVert}{\lVert h\rVert}
 \longrightarrow0.
\end{equation}
Equation \eqref{eq:state-Frechet-remainder} shows that $h\mapsto z[h]$ is the Fr\'echet derivative of
$\theta\mapsto x^\theta(\cdot,y)$ at $\bar\theta$.  Taking $h=\varepsilon\eta$ gives \eqref{eq:state-variation}.  Estimate
\eqref{eq:state-variation-bound} follows from
\eqref{eq:linearized-state-bound-h} with $h=\eta$.  Uniqueness follows by applying Gronwall's inequality to the difference of two solutions of \eqref{eq:state-variation}.

Finally, $z_y(S)=0$ because the same initial coordinate $\xi$ is used for all nearby policy parameters.  This argument is realization-wise and nowhere uses independence of $\omega$ and $\xi$.
\end{resultproof}

\subsection{Adjoint-based policy gradient}

\begin{resultproof}{Theorem~\ref{thm:policy-gradient}}
We divide the proof into four steps.

\smallskip
\noindent\emph{Step 1: closed-loop derivatives and the adjoint equation.}
Along the reference trajectory, define
\begin{align*}
 A_y(t)&:=D_xF_{\bar\theta}(t,x^{\bar\theta}(t,y),\omega),\\
 q_y(t)&:=\nabla_x\ell_{\bar\theta}(t,x^{\bar\theta}(t,y),\omega),\\
 C_y(t)&:=D_\theta F_{\bar\theta}(t,x^{\bar\theta}(t,y),\omega),\\
 r_y(t)&:=\nabla_\theta\ell_{\bar\theta}(t,x^{\bar\theta}(t,y),\omega).
\end{align*}
Thus $A_y$, $C_y$, $q_y$, and $r_y$ take values in $\R^{n\times n}$, $\R^{n\times q}$, $\R^n$, and $\R^q$, respectively.  The chain rule gives
\begin{align}
 A_y&=D_xf+D_ufD_x\pi_{\bar\theta},
 &C_y&=D_ufD_\theta\pi_{\bar\theta},\label{eq:closed-loop-chain-F}\\
 q_y&=\nabla_xL+(D_x\pi_{\bar\theta})^\top\nabla_uL,
 &r_y&=(D_\theta\pi_{\bar\theta})^\top\nabla_uL.
 \label{eq:closed-loop-chain-cost}
\end{align}
Consequently, \eqref{eq:closed-loop-adjoint} is equivalently
\begin{equation}\label{eq:abstract-adjoint}
 -\dot p_y(t)=A_y(t)^\top p_y(t)-q_y(t),
 \qquad
 p_y(T)=-\nabla_xg(x^{\bar\theta}(T,y),\omega).
\end{equation}
For almost every $y$, the coefficients $A_y$ and $q_y$ are integrable in time.  After the change of variable $s=T-t$, equation \eqref{eq:abstract-adjoint} is a forward linear Carath\'eodory equation.  Standard linear ODE theory for integrable coefficients \citep{coddington1955theory} gives a unique absolutely continuous solution of its integral equation.  Picard iteration of that integral equation consists only of integrals and products of jointly measurable functions, so the iterates are jointly measurable in $(t,y)$; their pointwise limit supplies a jointly measurable version of $p_y$.

Taking norms in the variation-of-constants formula and using
$\lVert A_y\rVert\le a$ and $\lVert q_y\rVert\le b$ gives
\begin{equation}\label{eq:adjoint-bound}
 \sup_{t\in[S,T]}\lVert p_y(t)\rVert
 \le
 \left(e(y)+\int_S^T b(s,y)\,ds\right)
 \exp\!\left(\int_S^T a(s,y)\,ds\right).
\end{equation}
This estimate also proves uniqueness: the difference of two solutions satisfies the homogeneous terminal-value equation and is zero by Gronwall's inequality.

\smallskip
\noindent\emph{Step 2: derivative of one realization cost.}
Fix $y$ as above and a direction $\eta\in\R^q$.  By Lemma~\ref{lem:state-sensitivity}, the state derivative $z_y$ satisfies
\begin{equation}\label{eq:proof-state-variation}
 \dot z_y=A_yz_y+C_y\eta,
 \qquad z_y(S)=0.
\end{equation}
Applying the ordinary chain rule to the terminal cost and the integral chain rule to the running cost along the state variation in \eqref{eq:proof-state-variation} gives
\begin{equation}\label{eq:realization-cost-derivative}
 DJ_y(\bar\theta)[\eta]
 =\nabla_xg(x^{\bar\theta}(T,y),\omega)^\top z_y(T)
 +\int_S^T q_y(t)^\top z_y(t)\,dt
 +\int_S^T r_y(t)^\top\eta\,dt.
\end{equation}
The differentiation is legitimate because the state remainder from Lemma~\ref{lem:state-sensitivity} is
$o(\lVert h\rVert)$ uniformly in time, while the derivative bounds in Assumption~\ref{ass:regularity} provide integrable time dominators for the mean-value remainders of $g$ and $\ell_\theta$.

\smallskip
\noindent\emph{Step 3: elimination of the state variation.}
The sign convention in \eqref{eq:abstract-adjoint} gives
\begin{equation}\label{eq:adjoint-proof-identities}
 q_y=A_y^\top p_y+\dot p_y,
 \qquad
 \nabla_xg(x^{\bar\theta}(T,y),\omega)=-p_y(T).
\end{equation}
Substitute \eqref{eq:adjoint-proof-identities} into the first two terms of
\eqref{eq:realization-cost-derivative}.  Since $p_y$ and $z_y$ are absolutely continuous, integration by parts is valid and yields
\begin{equation}\label{eq:adjoint-integration-by-parts}
\begin{aligned}
 &\nabla_xg(x^{\bar\theta}(T,y),\omega)^\top z_y(T)
 +\int_S^T q_y^\top z_y\,dt\\
 &\quad=-p_y(T)^\top z_y(T)
 +\int_S^T\left(p_y^\top A_yz_y+\dot p_y^\top z_y\right)dt\\
 &\quad=-p_y(T)^\top z_y(T)
 +\bigl[p_y^\top z_y\bigr]_{S}^{T}
 +\int_S^T p_y^\top(A_yz_y-\dot z_y)\,dt\\
 &\quad=-p_y(S)^\top z_y(S)-\int_S^T p_y^\top C_y\eta\,dt\\
 &\quad=-\int_S^T p_y^\top C_y\eta\,dt.
\end{aligned}
\end{equation}
In \eqref{eq:adjoint-integration-by-parts}, the terminal terms cancel by the terminal condition in \eqref{eq:closed-loop-adjoint}, and the initial term vanishes because $z_y(S)=0$.  Combining this identity with
\eqref{eq:realization-cost-derivative}, \eqref{eq:closed-loop-chain-F}, and
\eqref{eq:closed-loop-chain-cost} gives
\begin{align}
 DJ_y(\bar\theta)[\eta]
 &=\int_S^T\left(r_y-C_y^\top p_y\right)^\top\eta\,dt\nonumber\\
 &=\int_S^T
 \eta^\top(D_\theta\pi_{\bar\theta})^\top
 \left(\nabla_uL-D_uf^\top p_y\right)dt.
 \label{eq:realization-policy-gradient}
\end{align}

\smallskip
\noindent\emph{Step 4: differentiation under the joint law.}
By \eqref{eq:adjoint-bound}, $\lVert C_y\rVert\le c$, and
$\lVert r_y\rVert\le d$, the absolute value of
\eqref{eq:realization-policy-gradient} is bounded by
\begin{equation}\label{eq:realization-gradient-envelope}
 \left[
 \left(e(y)+\int_S^T b(t,y)\,dt\right)
 \exp\!\left(\int_S^T a(t,y)\,dt\right)
 \int_S^T c(t,y)\,dt
 +\int_S^T d(t,y)\,dt
 \right]\lVert\eta\rVert.
\end{equation}
The bracket in \eqref{eq:realization-gradient-envelope} is precisely the integrable envelope in
\eqref{eq:dominating-envelope}.  To control the full cost difference, use the mean-value formula separately for the terminal cost, the state argument of the running cost, and the policy-parameter argument of the running cost.  Write
$\Delta_hJ_y:=J_y(\bar\theta+h)-J_y(\bar\theta)$.  Together with \eqref{eq:state-lipschitz-bound}, this gives, for all sufficiently small nonzero $h$,
\begin{align}
 \frac{|\Delta_hJ_y|}{\lVert h\rVert}
 &\le
 \left(e(y)+\int_S^T b(t,y)\,dt\right)
 \frac{\sup_t\lVert x^{\bar\theta+h}(t,y)-x^{\bar\theta}(t,y)\rVert}{\lVert h\rVert}
 +\int_S^T d(t,y)\,dt\\
 &\le
 \left(e(y)+\int_S^T b(t,y)\,dt\right)
 \exp\!\left(\int_S^T a(t,y)\,dt\right)
 \int_S^T c(t,y)\,dt
 +\int_S^T d(t,y)\,dt.
 \label{eq:cost-difference-envelope}
\end{align}
The final expression in \eqref{eq:cost-difference-envelope} is the integrable envelope in
\eqref{eq:dominating-envelope}.  For almost every $y$, define the normalized realization remainder by
\begin{equation}\label{eq:realization-cost-remainder}
 \frac{J_y(\bar\theta+h)-J_y(\bar\theta)-DJ_y(\bar\theta)[h]}{\lVert h\rVert}.
\end{equation}
Step 2 and Lemma~\ref{lem:state-sensitivity} show that the expression in \eqref{eq:realization-cost-remainder} converges to zero.  Its absolute value is bounded by a fixed multiple of the same integrable envelope.  The dominated convergence theorem therefore permits differentiation under the $\nu$-integral.  The bound on \eqref{eq:realization-policy-gradient} also permits Fubini's theorem for the time--realization integral; see \citet{folland1999real} for both results.  Hence
\begin{equation}\label{eq:average-realization-derivative}
 D\cJ_\nu(\bar\theta)[\eta]
 =\int_\Y DJ_y(\bar\theta)[\eta] \,\nu(dy).
\end{equation}
Substituting \eqref{eq:realization-policy-gradient} into \eqref{eq:average-realization-derivative} proves
\eqref{eq:policy-gradient}.  Because the bound is linear in $\lVert\eta\rVert$ and the realization remainders are dominated uniformly for $h\to0$, the resulting derivative is Fr\'echet, not merely directional.  Finally, $\nabla_uH=D_uf^\top p_y-\nabla_uL$ gives the equivalent Hamiltonian form stated after \eqref{eq:policy-gradient}.
\end{resultproof}

\subsection{Policy-parameter stationarity}

\begin{resultproof}{Corollary~\ref{cor:projected-stationarity}}
First suppose that $\theta^*\in\interior\Theta$.  For every direction
$\eta\in\R^q$, both $\theta^*+\varepsilon\eta$ and
$\theta^*-\varepsilon\eta$ are feasible when $\varepsilon>0$ is sufficiently small.  Local minimality therefore gives
\begin{equation*}
 D\cJ_\nu(\theta^*)[\eta]\ge0
 \quad\text{and}\quad
 D\cJ_\nu(\theta^*)[-\eta]\ge0.
\end{equation*}
Because the derivative is linear in the direction, the second inequality is
$-D\cJ_\nu(\theta^*)[\eta]\ge0$.  Hence
$D\cJ_\nu(\theta^*)[\eta]=0$ for every $\eta$, which is equivalent to
$\nabla_\theta\cJ_\nu(\theta^*)=0$.  Using the Hamiltonian form of \eqref{eq:policy-gradient} and multiplying by $-1$ gives \eqref{eq:param-stationarity}.

Now suppose that $\Theta$ is closed and convex and that $\theta^*$ is a local minimizer relative to $\Theta$.  Fix any $\theta\in\Theta$ and set
$\eta:=\theta-\theta^*$.  Convexity implies
$\theta^*+s\eta=(1-s)\theta^*+s\theta\in\Theta$ for every
$s\in[0,1]$.  For all sufficiently small $s>0$, local minimality gives
\begin{equation*}
 \cJ_\nu(\theta^*+s\eta)-\cJ_\nu(\theta^*)\ge0.
\end{equation*}
Divide by $s$ and let $s\downarrow0$.  Differentiability yields
$D\cJ_\nu(\theta^*)[\eta]\ge0$, that is,
\begin{equation*}
 \left\langle\nabla_\theta\cJ_\nu(\theta^*),
 \theta-\theta^*\right\rangle\ge0.
\end{equation*}
This is \eqref{eq:param-VI}.

For a nonempty closed convex set, the Euclidean projection characterization gives, for $\vartheta\in\Theta$ and $d\in\R^q$,
\begin{equation*}
 \vartheta=\Pi_\Theta(\vartheta-d)
 \quad\Longleftrightarrow\quad
 \ip{d}{\theta-\vartheta}\ge0
 \qquad(\theta\in\Theta).
\end{equation*}
For the finite-law specialization~\eqref{eq:design-rule}, apply this identity with
$\vartheta=\theta^*$ and $d=\nabla_\theta\cJ_\nu(\theta^*)$.  The variational inequality for $\cJ_\nu$ is then equivalent to
$r_\Theta(\theta^*)=0$ in \eqref{eq:policy-stationarity-residual}.
\end{resultproof}

\subsection{Conditional moment propagation and cost}

We first justify the link between Assumptions~\ref{ass:regularity} and \ref{ass:lq-integrability}.  Here $A_i,Q_i,G_i\in\R^{n\times n}$, $B_i\in\R^{n\times m}$, and $R_i\in\R^{m\times m}$; the time-dependent matrices are understood pointwise almost everywhere.  On a bounded neighborhood of $K$, variation of constants gives $\sup_t\norm{x_i(t,\xi)}\le C_i\norm{\xi}$; the envelopes in Assumption~\ref{ass:regularity} may therefore be chosen with $a=O(1)$, $b,c,e=O(\norm{\xi})$, and $d=O(\norm{\xi}^2)$.  Thus \eqref{eq:dominating-envelope} is $O(\norm{\xi}^2)$ and is integrable by Assumption~\ref{ass:lq-integrability}, while the remaining conditions follow from the linear--quadratic data and the finite model marginal.

\begin{resultproof}{Proposition~\ref{prop:conditional-moments}}
Let $\Phi_{i,K}(t,s)$ be the state-transition matrix of
$\dot x=A_{i,K}(t)x$.  Essential boundedness of $A_{i,K}$ on the finite interval implies that $\Phi_{i,K}$ exists, is absolutely continuous in $t$, and is uniformly bounded on $[S,T]^2$ by standard linear-system theory \citep{sontag1998mathematical}.  For every initial state $\xi$,
\begin{equation}\label{eq:state-transition-solution}
 x_i(t,\xi)=\Phi_{i,K}(t,S)\xi.
\end{equation}
Substitution of \eqref{eq:state-transition-solution} into \eqref{eq:state-moment-definition} gives
\begin{equation}\label{eq:moment-transition-representation}
 M_i(t)
 =\Phi_{i,K}(t,S)M_i^0\Phi_{i,K}(t,S)^\top.
\end{equation}
The right-hand side is finite because $M_i^0$ is finite.  Differentiating
\eqref{eq:moment-transition-representation} for almost every $t$ and using
$\partial_t\Phi_{i,K}(t,S)=A_{i,K}(t)\Phi_{i,K}(t,S)$ gives
\begin{equation}\label{eq:moment-differentiated}
 \dot M_i=A_{i,K}M_i+M_iA_{i,K}^\top,
 \qquad M_i(S)=M_i^0.
\end{equation}
Equation \eqref{eq:moment-differentiated} is \eqref{eq:moment-forward}.  Equivalently, one may differentiate
$x_ix_i^\top$ before integration.  That interchange is justified by
$\lVert x_i(t,\xi)\rVert\le C_i\lVert\xi\rVert$ and
$\lVert\frac{d}{dt}(x_ix_i^\top)\rVert\le C_i'\lVert\xi\rVert^2$, together with the finite second moment of $\mu_i$.

For any symmetric matrix $C$ and any square-integrable random vector $X$,
\begin{equation}\label{eq:quadratic-trace-identity}
 \int X^\top C X\,d\mu
 =\int\tr(CXX^\top)\,d\mu
 =\tr\!\left(C\int XX^\top d\mu\right).
\end{equation}
Apply \eqref{eq:quadratic-trace-identity} at time $t$ with
$C=Q_i+K^\top R_iK$ and at the terminal time with $C=G_i$.  Since
$u_i=-Kx_i$,
\begin{equation}\label{eq:closed-loop-quadratic-cost}
 x_i^\top Q_ix_i+u_i^\top R_iu_i
 =x_i^\top(Q_i+K^\top R_iK)x_i.
\end{equation}
Tonelli's theorem \citep{folland1999real} applies to these nonnegative quadratic costs; equivalently, bounded coefficients and the finite second moment give absolute integrability.  Applying \eqref{eq:quadratic-trace-identity} to \eqref{eq:closed-loop-quadratic-cost}, integrating first over $\mu_i$, then over time, and finally summing with weights $\alpha_i$ produces the displayed trace formula for $\cJ_\nu(K)$.

No first moment appears because both the homogeneous dynamics and the cost are purely quadratic: all dependence on $\xi$ is through $\xi\xi^\top$.
\end{resultproof}

\subsection{Lyapunov evaluation and value}

\begin{resultproof}{Proposition~\ref{prop:policy-evaluation}}
Write
\begin{equation*}
 Q_{i,K}:=Q_i+K^\top R_iK.
\end{equation*}
Because $A_{i,K}$ and $Q_{i,K}$ are essentially bounded, the backward linear matrix equation \eqref{eq:P-Lyapunov} has a unique absolutely continuous solution.  Since $G_i$ and $Q_{i,K}$ are symmetric, transposing the equation shows that $P_i^\top$ solves the same terminal-value problem as $P_i$.  Uniqueness therefore implies $P_i(t)=P_i(t)^\top$ for all $t$.

Define $\widehat p_i:=-P_ix_i$.  Using
$\dot x_i=A_{i,K}x_i$ gives
\begin{equation}\label{eq:policy-evaluation-substitution}
\begin{aligned}
 -\dot{\widehat p}_i
 &=\dot P_i x_i+P_i\dot x_i\\
 &=\left(-A_{i,K}^\top P_i-P_iA_{i,K}-Q_{i,K}\right)x_i
   +P_iA_{i,K}x_i\\
 &=-A_{i,K}^\top P_ix_i-Q_{i,K}x_i\\
 &=A_{i,K}^\top\widehat p_i-Q_{i,K}x_i.
\end{aligned}
\end{equation}
At $t=T$,
$\widehat p_i(T)=-P_i(T)x_i(T)=-G_ix_i(T)$.  Hence
$-\widehat p_i(T)=G_ix_i(T)$.  Thus \eqref{eq:policy-evaluation-substitution} and the terminal identity show that $\widehat p_i$ satisfies the specialization of \eqref{eq:closed-loop-adjoint} to $y=(\omega_i,\xi)$.  Uniqueness of that linear terminal-value problem proves
$p_i=\widehat p_i=-P_ix_i$.

For the value identity, use the product rule for absolutely continuous matrix functions and substitute \eqref{eq:moment-forward} and \eqref{eq:P-Lyapunov}:
\begin{equation}\label{eq:trace-product-expansion}
\begin{aligned}
 \frac{d}{dt}\tr(P_iM_i)
 &=\tr(\dot P_iM_i)+\tr(P_i\dot M_i)\\
 &=\tr\!\left[
 (-A_{i,K}^\top P_i-P_iA_{i,K}-Q_{i,K})M_i
 \right]\\
 &\quad+\tr\!\left[P_i(A_{i,K}M_i+M_iA_{i,K}^\top)\right].
\end{aligned}
\end{equation}
The two terms containing $P_iA_{i,K}M_i$ cancel directly.  The remaining pair also cancels by cyclic invariance of the trace:
\begin{equation}\label{eq:trace-cyclic-cancellation}
 \tr(A_{i,K}^\top P_iM_i)
 =\tr(P_iM_iA_{i,K}^\top).
\end{equation}
Applying \eqref{eq:trace-cyclic-cancellation} to \eqref{eq:trace-product-expansion} gives
\begin{equation}\label{eq:trace-product-derivative}
 \frac{d}{dt}\frac12\tr(P_iM_i)
 =-\frac12\tr(Q_{i,K}M_i)
 \quad\text{for almost every }t.
\end{equation}
Integrating \eqref{eq:trace-product-derivative} from $S$ to $T$ gives
\begin{equation}\label{eq:integrated-trace-identity}
 \frac12\tr(P_i(S)M_i^0)
 =\frac12\tr(P_i(T)M_i(T))
 +\frac12\int_S^T\tr(Q_{i,K}M_i)\,dt.
\end{equation}
Substitute $P_i(T)=G_i$ into \eqref{eq:integrated-trace-identity}, multiply by $\alpha_i$, and sum over $i$.  The right-hand side is exactly the previously derived trace expression for $\cJ_\nu(K)$, proving \eqref{eq:trace-value-identity}.
\end{resultproof}

\subsection{Moment-based shared-gain gradient}

\begin{resultproof}{Theorem~\ref{thm:shared-gain-derivative}}
Fix a matrix variation $\Delta K\in\R^{m\times n}$ and define
$K_\varepsilon:=K+\varepsilon\Delta K$.  For sufficiently small $\varepsilon$, $K_\varepsilon$ defines a well-posed linear closed loop on the finite horizon.  Let
$x_i^\varepsilon(t,\xi)$ be the corresponding state and set
\begin{equation*}
 z_i(t,\xi)
 :=\left.\frac{d}{d\varepsilon}x_i^\varepsilon(t,\xi)\right|_{\varepsilon=0}.
\end{equation*}
Differentiating
$\dot x_i^\varepsilon=(A_i-B_iK_\varepsilon)x_i^\varepsilon$ gives
\begin{equation}\label{eq:LQ-gain-state-variation}
 \dot z_i=A_{i,K}z_i-B_i\Delta K\,x_i,
 \qquad z_i(S,\xi)=0.
\end{equation}
Differentiation of the linear integral equation shows that this derivative exists uniformly in time for each $\xi$.  Because the system is linear, both $x_i$ and $z_i$ are bounded by a constant times $\lVert\xi\rVert$; their quadratic products are therefore integrable under the finite-second-moment law $\mu_i$.

Let
$Q_{i,K}:=Q_i+K^\top R_iK$.  For a fixed initial state, differentiation of the quadratic cost gives
\begin{align}
 DJ_{(\omega_i,\xi)}(K)[\Delta K]
 &=x_i(T,\xi)^\top G_i z_i(T,\xi)
   +\int_S^T x_i^\top Q_{i,K}z_i\,dt\nonumber\\
 &\quad+\int_S^T x_i^\top K^\top R_i\Delta K\,x_i\,dt.
 \label{eq:LQ-cost-direct-variation}
\end{align}
The last term is the direct derivative of
$\tfrac12x_i^\top K^\top R_iKx_i$ with the state held fixed; symmetry of $R_i$ combines the two product-rule terms.

Let $p_i=-P_ix_i$, as established in Proposition~\ref{prop:policy-evaluation}.  Specializing \eqref{eq:closed-loop-adjoint} to the linear--quadratic data gives
\begin{equation}\label{eq:LQ-adjoint-identities}
 Q_{i,K}x_i=A_{i,K}^\top p_i+\dot p_i,
 \qquad G_ix_i(T)=-p_i(T).
\end{equation}
Use \eqref{eq:LQ-adjoint-identities} in the first two terms of
\eqref{eq:LQ-cost-direct-variation}.  Integration by parts and
$z_i(S)=0$ give
\begin{equation}\label{eq:LQ-adjoint-integration-by-parts}
\begin{aligned}
 &x_i(T)^\top G_i z_i(T)
 +\int_S^T x_i^\top Q_{i,K}z_i\,dt\\
 &\quad=-p_i(T)^\top z_i(T)
 +\int_S^T\left(p_i^\top A_{i,K}z_i+\dot p_i^\top z_i\right)dt\\
 &\quad=\int_S^T p_i^\top(A_{i,K}z_i-\dot z_i)\,dt\\
 &\quad=\int_S^T p_i^\top B_i\Delta K\,x_i\,dt.
\end{aligned}
\end{equation}
The last equality in \eqref{eq:LQ-adjoint-integration-by-parts} follows from
\eqref{eq:LQ-gain-state-variation}.  Since $p_i=-P_ix_i$, adding the direct term in
\eqref{eq:LQ-cost-direct-variation} yields
\begin{equation}\label{eq:LQ-realization-gain-derivative}
\begin{aligned}
 DJ_{(\omega_i,\xi)}(K)[\Delta K]
 &=\int_S^T
 x_i^\top\left(K^\top R_i-P_iB_i\right)\Delta K\,x_i\,dt\\
 &=\int_S^T
 \left\langle(R_iK-B_i^\top P_i)x_ix_i^\top,
 \Delta K\right\rangle_Fdt.
\end{aligned}
\end{equation}
The second equality in \eqref{eq:LQ-realization-gain-derivative} follows from
$\langle A,\Delta K\rangle_F=\tr(A^\top\Delta K)$ and cyclic invariance of the trace.

The integrand is bounded by a time-integrable constant times
$\lVert\xi\rVert^2\lVert\Delta K\rVert_F$.  Hence Fubini's theorem and differentiation under $\mu_i$ are justified by the finite second moment and dominated convergence \citep{folland1999real}.  Using
$\int x_ix_i^\top\,d\mu_i=M_i$ in \eqref{eq:LQ-realization-gain-derivative} gives
\begin{equation}\label{eq:LQ-averaged-gain-derivative}
\begin{aligned}
 D\cJ_\nu(K)[\Delta K]
 &=\sum_{i=1}^N\alpha_i\int_S^T
 \left\langle
 (R_iK-B_i^\top P_i)M_i,\Delta K
 \right\rangle_Fdt\\
 &=\left\langle
 \int_S^T\sum_{i=1}^N\alpha_i
 (R_iK-B_i^\top P_i)M_i\,dt,
 \Delta K
 \right\rangle_F.
\end{aligned}
\end{equation}
To obtain Fr\'echet differentiability, let $\Delta K\to0$ in $\R^{m\times n}$ and let $z_i[\Delta K]$ denote the corresponding solution of \eqref{eq:LQ-gain-state-variation}.  Variation of constants and Gronwall's inequality, uniformly for $K+\Delta K$ in a bounded neighborhood of $K$, give
\[
 \sup_{t\in[S,T]}
 \norm{x_i^{K+\Delta K}(t,\xi)-x_i^K(t,\xi)-z_i[\Delta K](t,\xi)}
 \le C_i\norm{\Delta K}_F^2\norm{\xi}.
\]
The corresponding cost remainder is bounded by
$C_i'\norm{\Delta K}_F^2\norm{\xi}^2$.  The finite second moment of $\mu_i$ therefore permits integration of the remainder and shows that it is $o(\norm{\Delta K}_F)$.  Hence \eqref{eq:LQ-averaged-gain-derivative} is the Fr\'echet derivative, and its integral coefficient is the Frobenius gradient in \eqref{eq:gain-derivative}.
\end{resultproof}

\subsection{Shared-gain stationarity}

\begin{resultproof}{Corollary~\ref{cor:shared-gain-stationarity}}
By Theorem~\ref{thm:shared-gain-derivative},
\begin{equation}\label{eq:constant-gain-derivative}
 D\cJ_\nu(K)[\Delta K]
 =\left\langle\nabla_K\cJ_\nu(K),\Delta K\right\rangle_F,
 \qquad
 \nabla_K\cJ_\nu(K)
 =\int_S^T\sum_{i=1}^N\alpha_i
 (R_iK-B_i^\top P_i)M_i\,dt.
\end{equation}
Because the gain space is $\R^{m\times n}$, both $\Delta K$ and $-\Delta K$ are feasible for every matrix direction.  Local minimality therefore gives
$D\cJ_\nu(K^*)[\Delta K]=0$ for every constant matrix $\Delta K$.  The only matrix whose Frobenius product with every $\Delta K$ vanishes is the zero matrix, so \eqref{eq:constant-gain-derivative} gives
\eqref{eq:shared-gain-stationarity}.
\end{resultproof}

\subsection{Shared-gain fixed-point identity under excitation}

\begin{resultproof}{Corollary~\ref{cor:shared-gain-fixed-point}}
Because $R_i(t)=R$ for every $i$ and $t$, the left-hand side of
\eqref{eq:shared-gain-stationarity} can be separated as
\begin{equation}\label{eq:shared-gain-stationarity-factorization}
\begin{aligned}
 &\int_S^T\sum_{i=1}^N\alpha_i
 (RK-B_i^\top P_i)M_i\,dt\\
 &\qquad=RK\left(\sum_{i=1}^N\alpha_i\int_S^TM_i(t)\,dt\right)
 -\sum_{i=1}^N\alpha_i\int_S^TB_i(t)^\top P_i(t)M_i(t)\,dt\\
 &\qquad=RK\mathsf S-\mathsf B.
\end{aligned}
\end{equation}
By \eqref{eq:shared-gain-stationarity-factorization}, stationarity is equivalent to
$RK\mathsf S=\mathsf B$.  If $\mathsf S$ is invertible, multiplication from the left by $R^{-1}$ and from the right by $\mathsf S^{-1}$ proves the equivalence with
$K=R^{-1}\mathsf B\mathsf S^{-1}$ in \eqref{eq:shared-gain-fixed-point}.
Moreover,
\begin{equation*}
 D(K)=-R^{-1}\nabla_K\cJ_\nu(K)\mathsf S(K)^{-1},
\end{equation*}
and hence, whenever the gradient is nonzero,
\begin{equation*}
 \ip{\nabla_K\cJ_\nu(K)}{D(K)}_F
 =-\norm{R^{-1/2}\nabla_K\cJ_\nu(K)\mathsf S(K)^{-1/2}}_F^2<0.
\end{equation*}
Thus \eqref{eq:anderson-moore-direction} is a strict descent direction.

By \eqref{eq:state-moment-definition}, each $M_i(t)$ is symmetric positive semidefinite.  Therefore $\mathsf S$ is also symmetric positive semidefinite.  For any
$v\in\R^n$, Tonelli's theorem and
\eqref{eq:state-moment-definition} give
\begin{equation}\label{eq:excitation-quadratic-form}
\begin{aligned}
 v^\top\mathsf Sv
 &=\sum_{i=1}^N\alpha_i\int_S^T v^\top M_i(t)v\,dt\\
 &=\sum_{i=1}^N\alpha_i\int_S^T
   \int_{\R^n}v^\top x_i(t,\xi)x_i(t,\xi)^\top v
   \,\mu_i(d\xi)dt\\
 &=\sum_{i=1}^N\alpha_i\int_S^T
   \int_{\R^n}\bigl(v^\top x_i(t,\xi)\bigr)^2
   \,\mu_i(d\xi)dt.
\end{aligned}
\end{equation}
A symmetric positive semidefinite matrix is invertible if and only if its quadratic form is strictly positive for every nonzero vector.  Identity \eqref{eq:excitation-quadratic-form} therefore proves the equivalence between
$\mathsf S\succ0$ and \eqref{eq:excitation-condition}.

Finally, suppose $M_{i_0}^0\succ0$ for some model index $i_0$ with
$\alpha_{i_0}>0$.  From \eqref{eq:moment-transition-representation},
\begin{equation}\label{eq:positive-moment-propagation}
 M_{i_0}(t)=\Phi_{i_0,K}(t,S)M_{i_0}^0\Phi_{i_0,K}(t,S)^\top.
\end{equation}
In \eqref{eq:positive-moment-propagation}, the state-transition matrix $\Phi_{i_0,K}(t,S)$ is invertible, with inverse
$\Phi_{i_0,K}(S,t)$.  Hence, for every $v\ne0$,
\begin{equation}\label{eq:positive-propagated-quadratic-form}
 v^\top M_{i_0}(t)v
 =\bigl(\Phi_{i_0,K}(t,S)^\top v\bigr)^\top
 M_{i_0}^0
 \bigl(\Phi_{i_0,K}(t,S)^\top v\bigr)>0.
\end{equation}
By \eqref{eq:positive-propagated-quadratic-form}, the $i_0$ contribution to $v^\top\mathsf Sv$ is therefore strictly positive, while all other contributions are nonnegative.  Thus
$\mathsf S\succ0$.  This condition is sufficient but not necessary because time integration and averaging across several singular moments can also produce a positive-definite $\mathsf S$.
\end{resultproof}

\raggedbottom
\section{Numerical implementation details and additional trajectory figures}\label{sec:numerical-appendix}

All computations were performed on a MacBook Pro with an Apple M4 Max chip and 128~GB of memory.  The implementation uses Python~3.11.16 with NumPy~2.4.2 and SciPy~1.17.0.  All state, adjoint, moment, and Lyapunov equations are integrated with {DOP853}.  The source code is available at \url{https://github.com/inody/averaged-cost-control}.

The nonlinear adjoint method uses a width-16, one-hidden-layer tanh policy with input $(t,x_1,x_2)$, scalar output, and 81 parameters.  A supervised fit to a time-dependent linear teacher followed by direct refinement at the mean initial state provides the warm start.  Joint-law refinement uses {L-BFGS-B}, the gradient in Theorem~\ref{thm:policy-gradient}, and coordinate bounds $[-8,8]^{81}$.  The iteration stops when the residual in \eqref{eq:policy-stationarity-residual} is at most $10^{-4}$ and the relative design-objective change is at most $10^{-7}$; the residual is then recomputed under tighter verification tolerances.  The largest verified residual among the four selected neural policies is $9.58\times10^{-5}$.  For averaged design, optimization and validation use separate $32\times8$ Cartesian empirical laws; nominal design uses the mean pair.  These empirical laws are held fixed throughout optimization and selection.  Among checkpoints satisfying this finite-design stationarity certificate and completing all validation simulations with finite costs, the checkpoint with the lowest validation cost is selected.  The selected policies are compared under a common final-evaluation law formed by 128 model-parameter nodes and 16 held-out initial-state particles.

The parameter bounds alone do not guarantee numerical admissibility.  If a design evaluation fails to produce a finite full-horizon state, cost, and adjoint inventory, the optimizer receives the finite penalty $10^{12}+\norm{\theta}^2$ and gradient $2\theta$.  Any validation failure rejects the checkpoint; costs are never averaged over only the successful trajectories.  The final neural and static-feedback evaluations use the same all-or-nothing rule on their complete finite inventories.

For both linearized design models, the common conditional initial moment is $M_i^0=M^0:=\bar\xi\bar\xi^\top+(r_\xi^2/3)I_2$ at every model node $i$, with $r_\xi=0.2$ for both oscillators.  Averaged gain design uses the same 32-node model-parameter rule as the nonlinear adjoint method.  Starting from $K_0=(1,4)$, the calculation applies Algorithm~\ref{alg:linear-design} with Armijo constant $10^{-4}$ and step contraction $1/2$, and stops when the residual in \eqref{eq:gain-stationarity-residual} is at most $10^{-8}$ under the verification tolerances.  The largest verified residual among the four reported gains is $9.27\times10^{-9}$.  To check the implementation of \eqref{eq:gain-derivative}, we compare its two components at $K_0=(1,4)$ with the centered finite differences $[\cJ_\nu(K_0+hE_j)-\cJ_\nu(K_0-hE_j)]/(2h)$, where $h=10^{-5}$ and $E_j$ is the $j$th coordinate direction.  The maximum absolute componentwise difference over the four linear design problems is $2.01\times10^{-9}$.  The full Anderson--Moore step is accepted throughout both Van der Pol calculations and the nominal Stuart--Landau calculation.  Averaged Stuart--Landau design uses step sizes $1/32$, $1/128$, and $1/4$ in its first three updates and full steps thereafter.  Nominal design uses $\bar\xi\bar\xi^\top$ during optimization.  Final costs for static linear feedback are computed from the nonlinear closed-loop systems using the same 128-by-16 laws as for time-dependent neural feedback.

The following displays use the same three-row state-and-control layout as Figure~\ref{fig:nonlinear-vdp-trajectories}.  Within each oscillator, the averaged and nominal panels use the same pairs of model parameters and initial states.

\begin{figure}[H]
\centering
\includegraphics[width=0.80\linewidth]{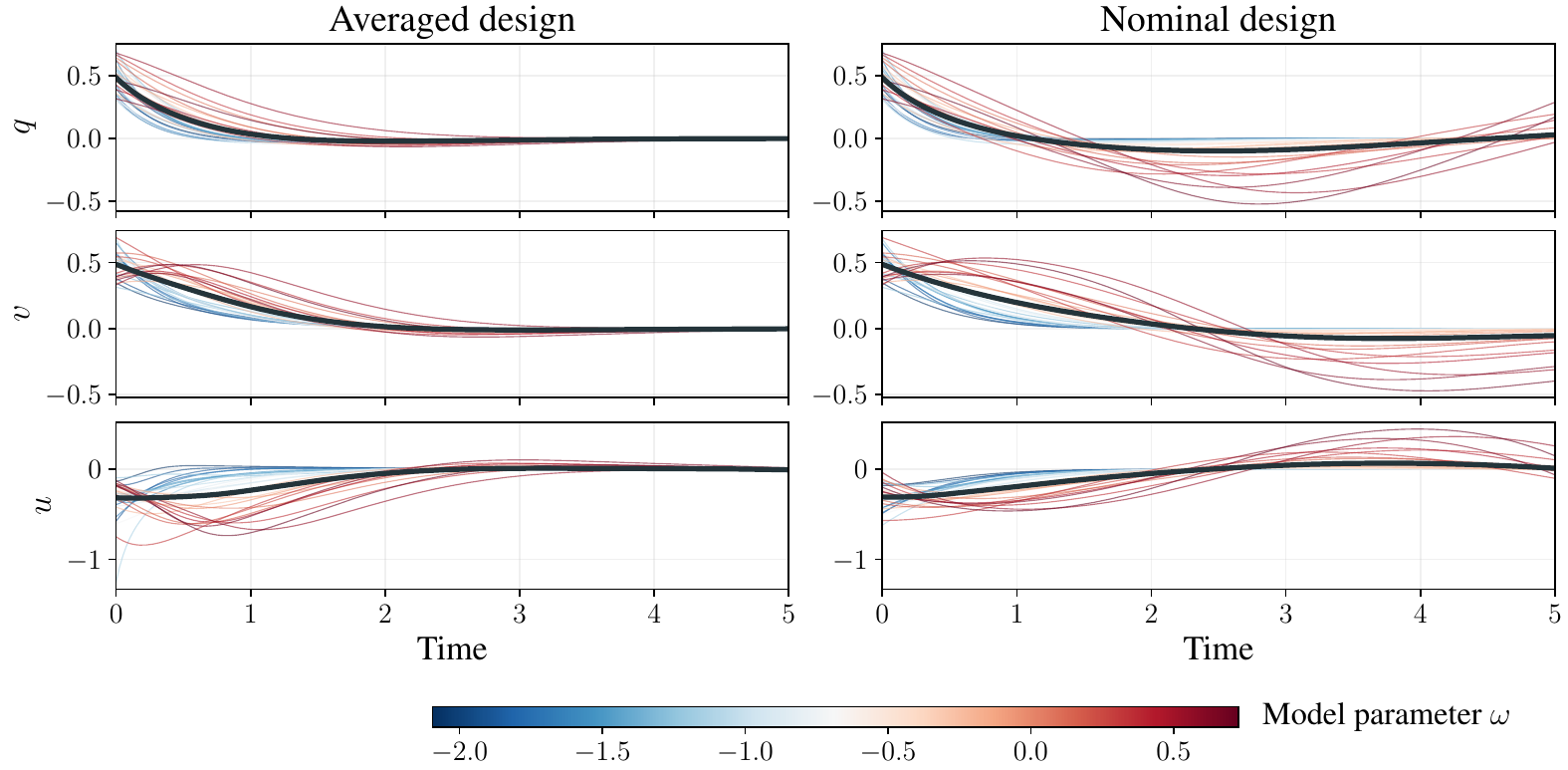}
\caption{Nonlinear Stuart--Landau trajectories under time-dependent neural feedback obtained by averaged design (left) and nominal design (right).  Both panels use the same 32 pairs of model parameters and initial states.  Rows show $q$, $v$, and $u$; color encodes the model parameter, and thick curves are equal-weight pointwise means.}
\label{fig:nonlinear-stuart-trajectories}
\end{figure}

\begin{figure}[H]
\centering
\includegraphics[width=0.80\linewidth]{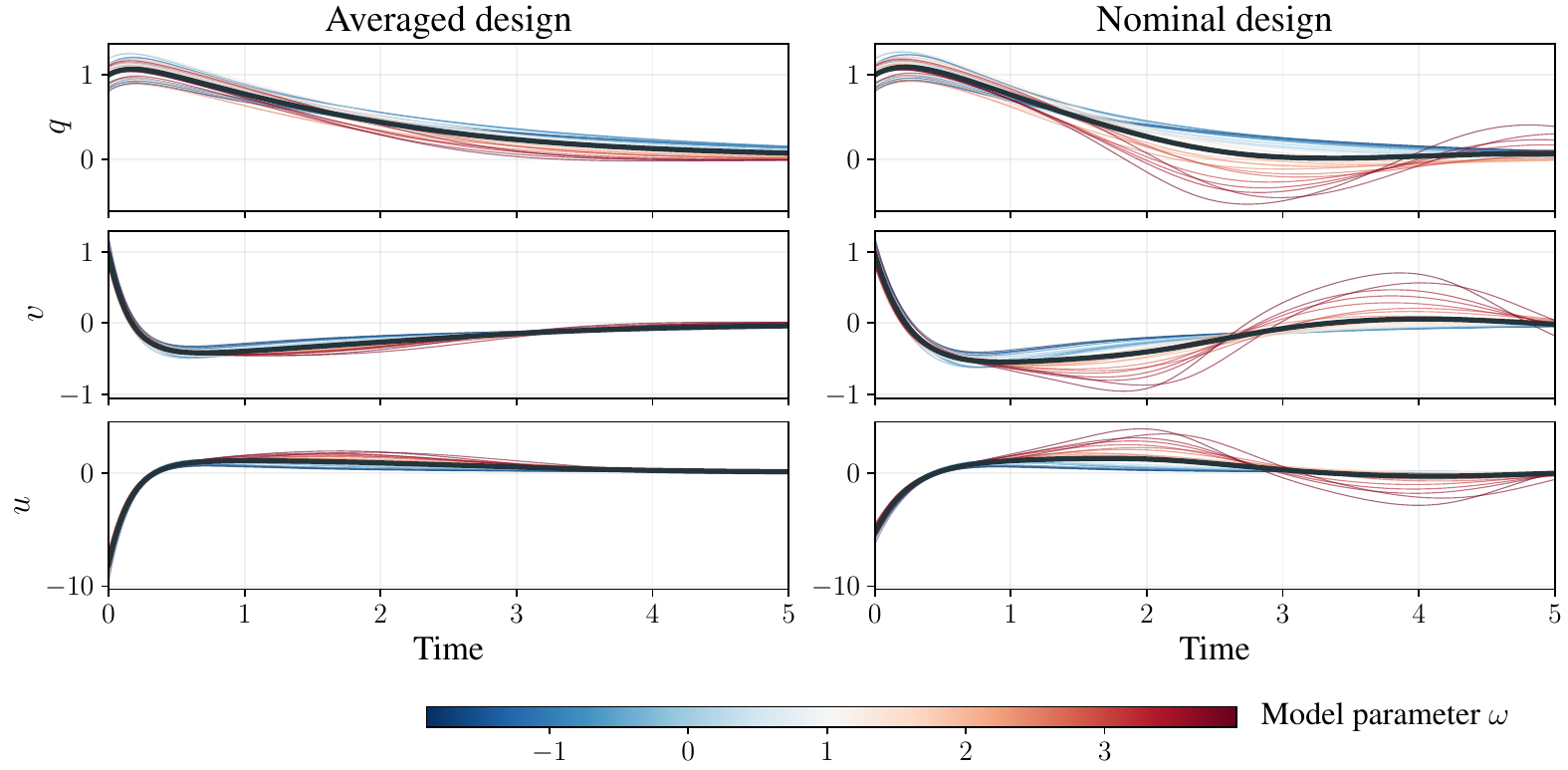}
\caption{Nonlinear Van der Pol trajectories under static linear feedback designed on the first-order linearization by averaged design (left) and nominal design (right).  Both panels use the same 32 pairs of model parameters and initial states.  Rows show $q$, $v$, and $u$; color encodes the model parameter, and thick curves are equal-weight pointwise means.  The gains are $K_{\mathrm{avg}}=(1.879565,\allowbreak 6.376037)$ and $K_{\mathrm{nom}}=(1.449509,\allowbreak 3.984057)$.}
\label{fig:nonlinear-vdp-static-trajectories}
\end{figure}

\begin{figure}[H]
\centering
\includegraphics[width=0.80\linewidth]{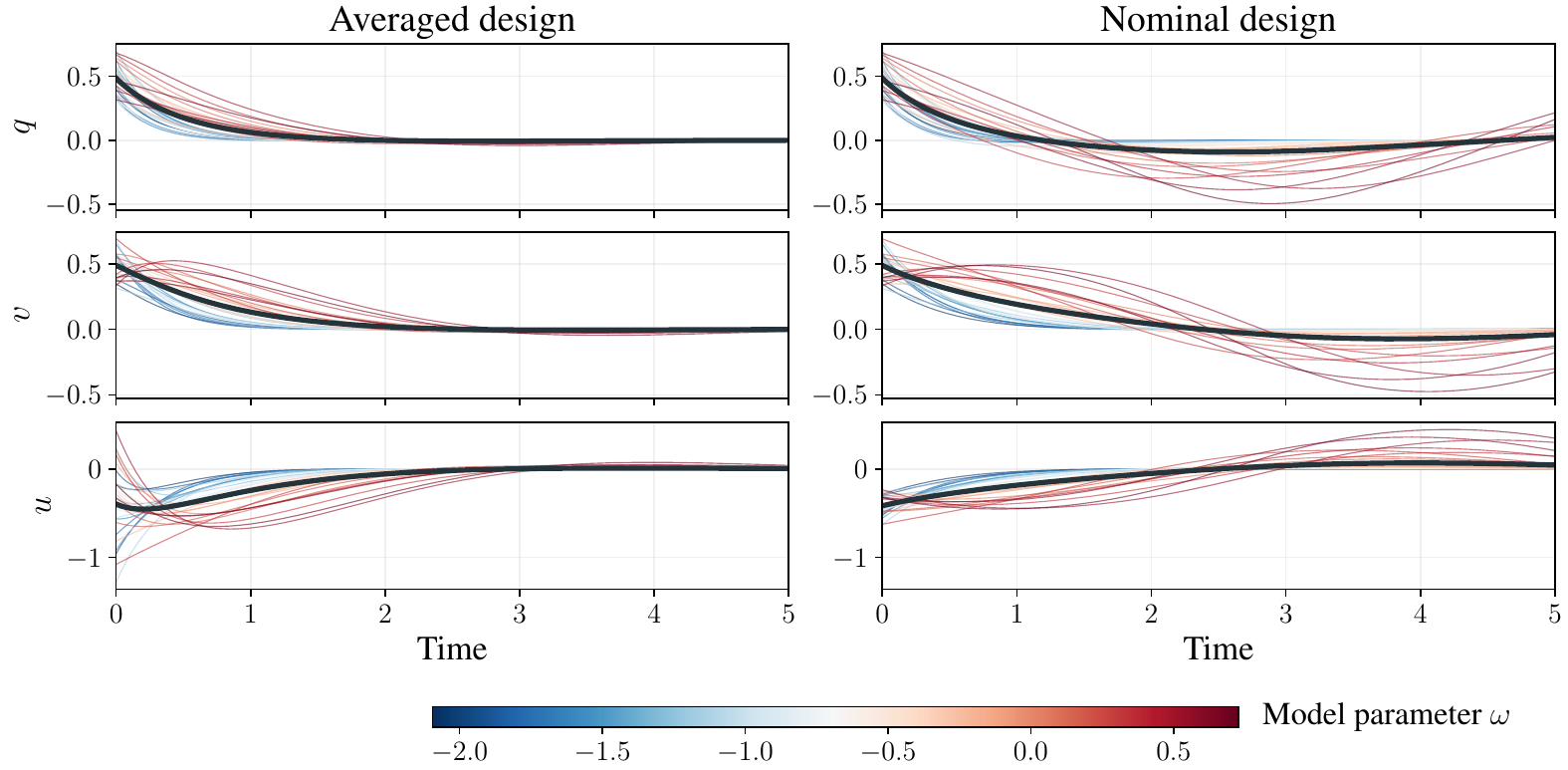}
\caption{Nonlinear Stuart--Landau trajectories under static linear feedback designed on the first-order linearization by averaged design (left) and nominal design (right).  Both panels use the same 32 pairs of model parameters and initial states.  Rows show $q$, $v$, and $u$; color encodes the model parameter, and thick curves are equal-weight pointwise means.  The gains are $K_{\mathrm{avg}}=(-1.994817,\allowbreak 2.801133)$ and $K_{\mathrm{nom}}=(-0.138222,\allowbreak 0.992785)$.}
\label{fig:nonlinear-stuart-static-trajectories}
\end{figure}

\begin{figure}[H]
\centering
\includegraphics[width=0.80\linewidth]{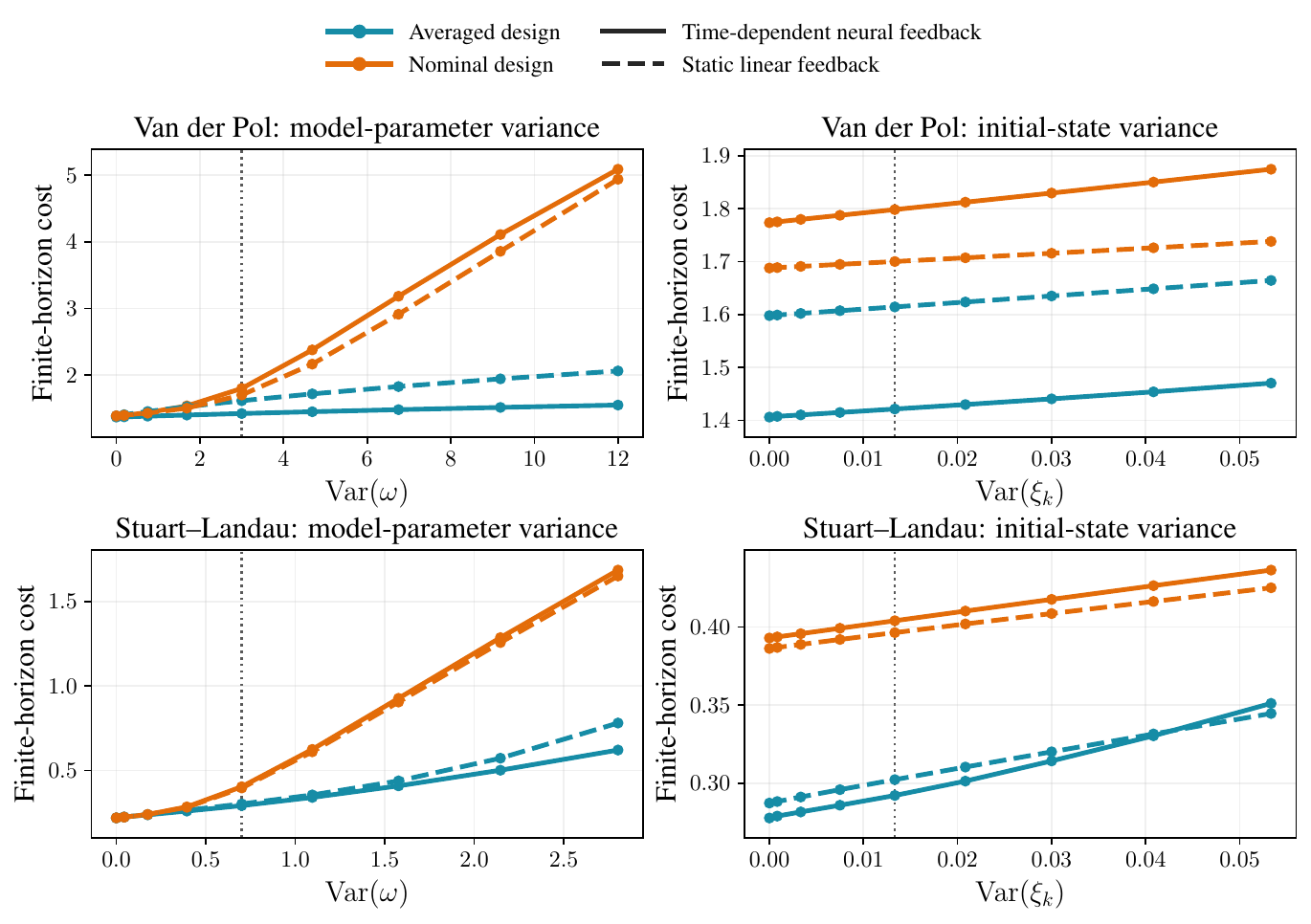}
\caption{Finite-horizon costs as the uncertainty variance changes.  At each point, averaged-design controllers are recomputed for the corresponding distribution, and both designs are evaluated under that distribution.  Rows identify the Van der Pol and Stuart--Landau systems; columns vary the model-parameter variance and component variance of the initial state while holding the other variance at the setting in Table~\ref{tab:evaluation-costs}.  Expectations use 128 model-parameter nodes for each nondegenerate model distribution and a symmetric $4\times4$ Gauss--Legendre initial-state rule.  Blue and orange denote averaged and nominal design; solid and dashed curves denote time-dependent neural feedback and static linear feedback.  Dotted lines mark the distributions used in Table~\ref{tab:evaluation-costs}.}
\label{fig:uncertainty-variance-sensitivity}
\end{figure}

Figure~\ref{fig:uncertainty-variance-sensitivity} reports the finite-horizon cost as the uncertainty variance changes.  At each plotted variance, we recompute the averaged-design neural policy and static linear gain using the corresponding model-parameter and initial-state distribution, and then evaluate these controllers and the nominal-design controllers under that same distribution.  The model-parameter variance and the component variance of the initial state are varied separately while preserving their means, and the other uncertainty source is held at the setting used in Table~\ref{tab:evaluation-costs}.  Because the mean pair does not change, the nominal-design controllers remain the same across both sweeps.

Increasing either variance raises all four costs.  Averaged design has lower cost throughout the initial-state variance sweeps for both feedback classes and both oscillators.  The model-parameter variance produces a stronger separation: the nominal costs increase more rapidly, and averaged design has lower cost at most settings.  The few reversals are confined to small model-parameter variances and remain small.  Thus the benefit of accounting for the uncertainty law is most pronounced as model-parameter variance increases, while the ordering is largely unchanged as initial-state variance varies.

\end{document}